# Exact Solutions to Cubic Duffing Equation

# by Leaf Functions Under Free Vibration

Kazunori Shinohara[1]

[1]Department of Mechanical System Engineering, Daido University
10-3 Takiharu-cho, Minami-ku, Nagoya 457-8530, Japan

*Corresponding author. Email: shinohara@06.alumni.u-tokyo.ac.jp

*Abstract*

Exact solutions with the initial conditions are presented in the cubic duffing equation. These exact solutions are expressed in terms of the leaf function and the trigonometric function. The leaf functions: $r=sleaf_n(t)$ or $r=cleaf_n(t)$ satisfy the ordinary differential equation: $dx^2/dt^2=-nr^{2n-1}$. The second order differential of the leaf function is equal to $-n$ times the function raised to the $(2n-1)$ power of the leaf function. By using the Leaf functions, the exact solutions of the cubic Duffing equation can be derived under several conditions. In this paper, seven types of the exact solutions are presented based on the leaf functions. This paper reveals the derivation of the seven exact solutions. Finally, the waves based on the exact solutions are visualized by the numerical results in the graphs

*1. Introduction*

The Duffing equation is the ordinary differential equation proposed by Georg Duffing[1][2]. The equation can be applied to various physical phenomenon (for e.g., the oscillator[3][4][5], pendulum[6], and circuit[7]). Methods based on numerical analysis were proposed to obtain the solution of the Duffing equation, and the calculation accuracy was discussed in previous works[3][8][9][10][11]. The author proposed the leaf function for this purpose[12][13]. Differential equations include ordinary differential equations (ODEs), which depend on one unknown variable, and partial differential equations, which depend on several unknown variables. In a previous study, the authors investigated an ODE in which the second order differential $dx^2/dt^2$ is equal to $-nx^{2n-1}$. The equation is as follows:

$$\frac{d^2x(t)}{dl^2} = -n \cdot x(t)^{2n-1} \qquad (1.1)$$

$$x(0) = 0 \qquad (1.2)$$

$$\frac{dx(0)}{dt} = 1 \qquad (1.3)$$

The graph of an ODE can be obtained numerically. In the graph, the horizontal and vertical axes represent variables *t* and *x*, respectively. Regular curves with a specific period and amplitude occur in all graphs. For *n = 1*, Eq. (1.1)–(1.3) represent the trigonometric function *sin*. For *n = 2*, these equations represent the lemniscatic elliptic function *sl*. However, the above equations have not yet been reported for *n = 3*. Therefore, the author named the function needed to obtain the curve groups as the leaf function $x(t)=sleaf_2(t)$, and continued studying this function[12][13]. The inverse leaf function is defined as follows:

$$arcsleaf_n(x) = \int_0^x \frac{1}{\sqrt{1-u^{2n}}} du = t \qquad (1.4)$$

The leaf function is defined as the function $cleaf_n(t)$. The ODEs and initial conditions that satisfy the function $x(t)=cleaf_2(t)$ are as follows:

$$\frac{d^2x(t)}{dl^2} = -n \cdot x(t)^{2n-1} \qquad (1.5)$$

$$x(0) = 1 \qquad (1.6)$$

$$\frac{dx(0)}{dt} = 0 \tag{1.7}$$

The inverse leaf function of the function $cleaf_n(t)$ is defined as follows:

$$arccleaf_n(x) = \int_x^1 \frac{1}{\sqrt{1-u^{2n}}} du = t \tag{1.8}$$

The periodicity of the functions $sleaf_n(t)$ and $cleaf_n(t)$ depends on the parameter $n$. The constant values of the periodicity are defined as follows:

$$2\pi_n = 4\int_0^1 \frac{1}{\sqrt{1-u^{2n}}} du \quad (n=1,2,3,\cdots) \tag{1.9}$$

The constant values $2\pi_n$ represent one periodicity with respect to arbitrary parameters $n$. The numerical results (for $n = 1, 2, 3$) are summarized in Table 1.

Table 1 Values of constant $\pi_n$

| $n$ | $\pi_n$ |
| --- | --- |
| 1 | $\pi_1 = 3.142\cdots$ |
| 2 | $\pi_2 = 2.622\cdots$ |
| 3 | $\pi_3 = 2.429\cdots$ |

The inverse leaf function for $n = 2$ is as follows:

$$arcsleaf_2(x) = \int_0^x \frac{1}{\sqrt{1-u^4}} du = t \tag{1.10}$$

$$arccleaf_2(x) = \int_x^1 \frac{1}{\sqrt{1-u^4}} du = t \tag{1.11}$$

Using Eq. (1.10) and Eq. (1.11), the numerical data between the parameters $x$ and $t$ can be obtained by numerical analyses. These data are also summarized in Table 1. The curves of the leaf functions $sleaf_2(t)$ and $cleaf_2(t)$ and the integral leaf functions $\int_0^t sleaf_2(u)du$ and $\int_0^t cleaf_2(u)du$ are shown in Fig. 1. The values of the integral functions can be obtained by numerical integration. The periodicity of the leaf functions $sleaf_2(t)$ and $cleaf_2(t)$ is $2\pi_2$. The mathematical description is as follows:

$$sleaf_2(t + 2\pi_2) = sleaf_2(t) \tag{1.12}$$

$$cleaf_2(t + 2\pi_2) = cleaf_2(t) \tag{1.13}$$

$$sleaf_2(m\pi_2) = 0 \quad (m = 0, \pm 1, \pm 2, \pm 3 \cdots) \tag{1.14}$$

$$sleaf_2\left(\frac{\pi_2}{2}(4m - 3)\right) = 1 \quad (m = 0, \pm 1, \pm 2, \pm 3 \cdots) \tag{1.15}$$

$$sleaf_2\left(\frac{\pi_2}{2}(4m - 1)\right) = -1 \quad (m = 0, \pm 1, \pm 2, \pm 3 \cdots) \tag{1.16}$$

$$cleaf_2\left(\frac{\pi_2}{2}(2m - 1)\right) = 0 \quad (m = 0, \pm 1, \pm 2, \pm 3 \cdots) \tag{1.17}$$

$$cleaf_2(2m\pi_2) = 1 \quad (m = 0, \pm 1, \pm 2, \pm 3 \cdots) \tag{1.18}$$

$$cleaf_2(\pi_2(2m - 1)) = -1 \quad (m = 0, \pm 1, \pm 2, \pm 3 \cdots) \tag{1.19}$$

In this paper, using the leaf functions: $sleaf_2(t)$ and $cleaf_2(t)$ for $n=2$, seven types of the exact solutions are presented for the cubic Duffing equation. In each case, the mathematical derivations and the numerical results of the seven types are shown in detail. Thereafter, the features of the waveforms are discussed.

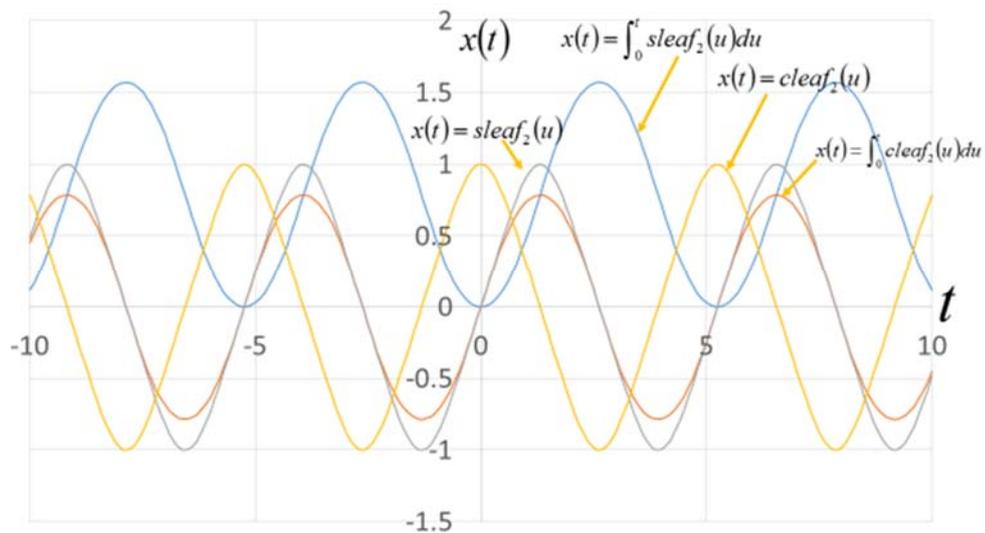

Fig. 1 Waves of $sleaf_2(t)$, $cleaf_2(t)$, $\int_0^t sleaf_2(u)du$ and $\int_0^t cleaf_2(u)du$

Table 2 Numerical data of *sleaf₂(t)*, *cleaf₂(t)*, $\int_0^t sleaf_2(u)du$ and $\int_0^t cleaf_2(u)du$

(All results have been rounded to no more than six significant figures)

| t | sleaf₂(t) | cleaf₂(t) | $\int_0^t sleaf_2(u)du$ | $\int_0^t cleaf_2(u)du$ |
|---|---|---|---|---|
| -10.0 | 0.48547 | 0.78647 | 0.11895 | 0.45195 |
| -9.0 | 0.96908 | -0.17718 | 0.96076 | 0.76969 |
| -8.0 | 0.13382 | -0.98225 | 1.56184 | 0.13303 |
| -7.0 | -0.81960 | -0.44312 | 1.20252 | -0.68658 |
| -6.0 | -0.73186 | 0.54991 | 0.28262 | -0.63179 |
| -5.0 | 0.24402 | 0.94212 | 0.02979 | 0.23935 |
| -4.0 | 0.99553 | 0.06691 | 0.71858 | 0.78315 |
| -3.0 | 0.37717 | -0.86655 | 1.49942 | 0.36067 |
| -2.0 | -0.61286 | -0.67373 | 1.37828 | -0.54982 |
| -1.0 | -0.90768 | 0.31073 | 0.48411 | -0.73704 |
| 0.0 | 0.00000 | 1.00000 | 0.00000 | 0.00000 |
| 1.0 | 0.90768 | 0.31073 | 0.48411 | 0.73704 |
| 2.0 | 0.61285 | -0.67373 | 1.37828 | 0.54982 |
| 3.0 | -0.37717 | -0.86655 | 1.49942 | -0.36067 |
| 4.0 | -0.99553 | 0.06691 | 0.71858 | -0.78316 |
| 5.0 | -0.24403 | 0.94212 | 0.02979 | -0.23935 |
| 6.0 | 0.73186 | 0.54991 | 0.28262 | 0.63179 |
| 7.0 | 0.81960 | -0.44312 | 1.20252 | 0.68657 |
| 8.0 | -0.13382 | -0.98225 | 1.56184 | -0.13303 |
| 9.0 | -0.96908 | -0.17718 | 0.96076 | -0.76970 |
| 10.0 | -0.48547 | 0.78647 | 0.11895 | -0.45196 |

## 2. Exact Solutions of Cubic Duffing Equation Using Leaf Functions

We try to apply the leaf function to the Duffing equation[12][13]. The Duffing equation is given as follows:

$$\frac{d^2x(t)}{dt^2} + \alpha x(t) + \beta x(t)^3 = 0 \tag{2.1}$$

For mechanical vibration, the above equation represents the free vibration by a nonlinear spring. The variable *x(t)* represents the unknown function and depends on the parameter *t*. Differential operators *dx(t)/dt* and *d²x(t)/dt²* represent the first and second order differential, respectively. The symbols α, β, and ω represent coefficients that do not depend on time *t*. In mechanical engineering fields, the above equation is regarded as the mathematical model for the nonlinear vibration. In the left side of Eq. (2.1), the first, second, and third terms represent inertia, stiffness, and nonlinear stiffness, respectively. By using the leaf functions, seven types of exact solutions can be set and then be derived the equation that the solution satisfies. In this paper, types (I)–(VII) are defined for the exact solutions, and the ordinary difference equations and the initial conditions are given as follows:

・Type (I)   (See Appendix I for details)

The exact solution: $x(t) = A\cos\left(\int_0^{\omega t+\phi} cleaf_2(u)du\right)$ (2.2)

The ordinary differential equation: $\dfrac{d^2x(t)}{dt^2} - 3\omega^2 x(t) + 4\left(\dfrac{\omega}{A}\right)^2 x(t)^3 = 0$ (2.3)

The initial position: $x(0) = A\cos\left(\int_0^{\phi} cleaf_2(t)dt\right)$ (2.4)

The initial velocity: $\dfrac{dx(0)}{dt} = -A\cdot\omega\cdot cleaf_2(\phi)\cdot\sin\left(\int_0^{\phi} cleaf_2(t)dt\right)$ (2.5)

・Type (II)   (See Appendix II for details)

The exact solution: $x(t) = A\sin\left(\int_0^{\omega t+\phi} cleaf_2(u)du\right)$ (2.6)

The ordinary differential equation: $\dfrac{d^2x(t)}{dt^2} + 3\omega^2 x(t) - 4\left(\dfrac{\omega}{A}\right)^2 x(t)^3 = 0$ (2.7)

The initial position: $x(0) = A\sin\left(\int_0^{\phi} cleaf_2(t)dt\right)$ (2.8)

The initial velocity: $\dfrac{dx(0)}{dt} = A\cdot\omega\cdot cleaf_2(\phi)\cdot\cos\left(\int_0^{\phi} cleaf_2(t)dt\right)$ (2.9)

・Type (III)   (See Appendix III for details)

The exact solution: $x(t) = A\cos\left(\int_0^{\omega t+\phi} sleaf_2(u)du\right) + A\sin\left(\int_0^{\omega t+\phi} sleaf_2(u)du\right)$ (2.10)

The ordinary differential equation: $\dfrac{d^2x(t)}{dt^2} - 3\omega^2 x(t) + 2\left(\dfrac{\omega}{A}\right)^2 x(t)^3 = 0$ (2.11)

The initial position: $x(0) = \sqrt{2}A\cos\left(\int_0^\phi sleaf_2(t)dt - \frac{\pi}{4}\right)$ (2.12)

The initial velocity: $\frac{dx(0)}{dt} = \sqrt{2}A \cdot \omega \cdot sleaf_2(\phi) \cdot \cos\left(\int_0^\phi sleaf_2(t)dt + \frac{\pi}{4}\right)$ (2.13)

· Type (IV)  (See Appendix IV for details)

The exact solution: $x(t) = A\cos\left(\int_0^{\omega t+\phi} sleaf_2(u)du\right) - A\sin\left(\int_0^{\omega t+\phi} sleaf_2(u)du\right)$ (2.14)

The ordinary differential equation: $\frac{d^2x(t)}{dt^2} + 3\omega^2 x(t) - 2\left(\frac{\omega}{A}\right)^2 x(t)^3 = 0$ (2.15)

The initial position: $x(0) = \sqrt{2}A\cos\left(\int_0^\phi sleaf_2(t)dt + \frac{\pi}{4}\right)$ (2.16)

The initial velocity: $\frac{dx(0)}{dt} = \sqrt{2}A \cdot \omega \cdot sleaf_2(\phi) \cdot \cos\left(\int_0^\phi sleaf_2(t)dt + \frac{\pi}{4}\right)$ (2.17)

· Type (V)  (See Appendix V for details)

The exact solution:

$x(t) = A\cos\left(\int_0^{\omega t+\phi} sleaf_2(u)du\right) + A\sin\left(\int_0^{\omega t+\phi} sleaf_2(u)du\right) + \sqrt{2}A\cos\left(\int_0^{\omega t+\phi} cleaf_2(u)du\right)$ (2.18)

The ordinary differential equation: $\frac{d^2x(t)}{dt^2} - 3\omega^2\left(1 + 2\sqrt{2}\right)x(t) + 2\frac{\omega^2}{A^2} x(t)^3 = 0$ (2.19)

The initial position: $x(0) = \sqrt{2}A\left\{\sin\left(\int_0^\phi sleaf_2(t)dt + \frac{\pi}{4}\right) + \cos\left(\int_0^\phi cleaf_2(t)dt\right)\right\}$ (2.20)

The initial velocity:

$\frac{dx(0)}{dt} = \sqrt{2}A \cdot \omega \cdot \left\{\cos\left(\int_0^\phi sleaf_2(t)dt + \frac{\pi}{4}\right) \cdot sleaf_2(\phi) - \sin\left(\int_0^\phi cleaf_2(t)dt\right) \cdot cleaf_2(\phi)\right\}$ (2.21)

· Type (VI)  (See Appendix VI for details)

The exact solution:

$x(t) = A\cos\left(\int_0^{\omega t+\phi} sleaf_2(u)du\right) + A\sin\left(\int_0^{\omega t+\phi} sleaf_2(u)du\right) - \sqrt{2}A\cos\left(\int_0^{\omega t+\phi} cleaf_2(u)du\right)$ (2.22)

The ordinary differential equation: $\frac{d^2x(t)}{dt^2} + 3\omega^2\left(2\sqrt{2} - 1\right)x(t) + 2\frac{\omega^2}{A^2} x(t)^3 = 0$ (2.23)

The initial position: $x(0) = \sqrt{2}A\left\{\sin\left(\int_0^\phi sleaf_2(t)dt + \frac{\pi}{4}\right) - \cos\left(\int_0^\phi cleaf_2(t)dt\right)\right\}$ (2.24)

The initial velocity:

$$\frac{dx(0)}{dt} = \sqrt{2} A \cdot \omega \cdot \left\{ \cos\left( \int_0^\phi sleaf_2(t)dt + \frac{\pi}{4} \right) \cdot sleaf_2(\phi) + \sin\left( \int_0^\phi cleaf_2(t)dt \right) \cdot cleaf_2(\phi) \right\} \quad (2.25)$$

・Type (VII)  (See Appendix VII for details)

The exact solution: $x(t) = A \cdot sleaf_2(\omega \cdot t + \phi) \cdot cleaf_2(\omega \cdot t + \phi)$ (2.26)

The ordinary differential equation: $\dfrac{d^2 x(t)}{dt^2} + 6\omega^2 x(t) - 2\left(\dfrac{\omega}{A}\right)^2 x(t)^3 = 0$ (2.27)

The initial position: $x(0) = A \cdot sleaf_2(\phi) \cdot cleaf_2(\phi)$ (2.28)

The initial velocity: $\dfrac{dx(0)}{dt} = A\omega \left\{ (cleaf_2(\phi))^2 - (sleaf_2(\phi))^2 \right\}$ (2.29)

The variables *x(t), A, t, ω, u,* and *Φ* represent displacement, amplitude, time, angular frequency, dummy variable and initial phase, respectively. The exact solution of the types (I)–(VII) is satisfied with the cubic duffing equation (Eq. (2.1)). These verifications are summarized in the appendix.

## 3. Numerical Results of Exact Solutions

### 3.1 Numerical Results of Exact Solution of Type(I)

In the case *A=1, ω=1, Φ=0* in Eq. (2.2), the wave of the exact solution of type (I) is compared to the wave of the function: $cleaf_2(t)$ and the wave of the function: $\int_0^t cleaf_2(u)du$. These waves are shown in Fig.2. The horizontal and vertical axes represent time and displacement *x(t)*, respectively. In the function: $cleaf_2(t)$, the amplitude and period become 1.0 and $2\pi_2$ (See the table 1), respectively. The center of the displacement becomes *x(t)=0*. In the function: $\int_0^t cleaf_2(u)du$, the amplitude becomes (See Appendix VIII):

$$\int_{0}^{\frac{\pi_2}{2}(2m-1)} cleaf_2(u)\,du = \pm\frac{\pi}{4}\,(\cong \pm 0.785398)\quad (m: integer) \tag{3.1}$$

The period becomes the constant $\pi_2$ (See the table 1). The center of the displacement becomes $x(t)=0$. Next, the wave obtained by type(I) exact solution is discussed. The minimum of the variable $x(t)$ is obtained as follows:

$$x\!\left(\frac{\pi_2}{2}(2m-1)\right) = \cos\!\left(\int_{0}^{\frac{\pi_2}{2}(2m-1)} cleaf_2(u)\,du\right) = \cos\!\left(\pm\frac{\pi}{4}\right) = \frac{1}{\sqrt{2}}\quad (m: integer) \tag{3.2}$$

The first order differential of the type( I ) solution is obtained as follows:

$$\frac{dx(t)}{dt} = -\sin\!\left(\int_{0}^{t} cleaf_2(u)\,du\right)\cdot cleaf_2(t) \tag{3.3}$$

By substituting $t = \dfrac{\pi_2}{2}(2m-1)$ into the Eq. (3.3), the Eq. (3.3) is satisfied with the following equation:

$$\begin{aligned}\frac{d}{dt}x\!\left(\frac{\pi_2}{2}(2m-1)\right) &= -\sin\!\left(\int_{0}^{\frac{\pi_2}{2}(2m-1)} cleaf_2(u)\,du\right)\cdot cleaf_2\!\left(\frac{\pi_2}{2}(2m-1)\right)\\ &= -\sin\!\left(\int_{0}^{\frac{\pi_2}{2}(2m-1)} cleaf_2(u)\,du\right)\cdot 0 = 0\end{aligned} \tag{3.4}$$

The Eq. (1.17) is applied to the above equation. Next, the maximum of variable $x(t)$ is obtained as follows:

$$x\!\left(\frac{\pi_2}{2}(2m)\right) = \cos\!\left(\int_{0}^{\frac{\pi_2}{2}(2m)} cleaf_2(u)\,du\right) = \cos(0) = 1.0\quad (m: integer) \tag{3.5}$$

By substituting $t = \dfrac{\pi_2}{2}(2m)$ into the Eq. (3.3), the Eq. (3.3) is satisfied with the following equation (See Appendix VIII):

$$\frac{d}{dt}x\left(\frac{\pi_2}{2}(2m)\right) = -\sin\left(\int_0^{\frac{\pi_2}{2}(2m)} cleaf_2(u)du\right) \cdot cleaf_2\left(\frac{\pi_2}{2}(2m)\right) = -\sin(0)\cdot(\pm 1) = 0 \quad (3.6)$$

In the type(I) solution, the range of the variable $x(t)$ becomes $\frac{1}{\sqrt{2}} \leq x(t) \leq 1$. The centers of the displacement and amplitude become $\frac{2+\sqrt{2}}{4}$ and $\frac{2-\sqrt{2}}{4}$, respectively.

In the solution of Type( I ), the integration function: $\int_0^{m\pi_2} cleaf_2(u)du$ represents the phase of the cosine function. As shown in Fig. 1, the range of the integration function: $\int_0^{m\pi_2} cleaf_2(u)du$ becomes the inequality: $-\frac{\pi}{4} \leq \int_0^{m\pi_2} cleaf_2(t)dt \leq \frac{\pi}{4}$ (See Appendix VIII). The periodicity of the integration function: $\int_0^{m\pi_2} cleaf_2(u)du$ becomes the constant $2\pi_2$ in the table 1. Therefore, in the case of $\omega=1$, $A=1$ and $\Phi=1$, the type ( I ) solution is as follows:

$$x(t) = \cos\left(\int_0^t cleaf_2(u)du\right) \quad (3.7)$$

The periodicity of the above solution becomes the constant $\pi_2$ because the cosine function is satisfied with the following equation:

$$\cos\left(\int_0^t cleaf_2(u)du\right) = \cos\left(-\int_0^t cleaf_2(u)du\right) \quad (3.8)$$

The ordinary difference equation is as follows:

$$\frac{d^2 x(t)}{dt^2} - 3x(t) + 4x(t)^3 = 0 \quad (3.9)$$

The second derivative $d^2x(t)/dt^2$ is as follows(See Apendix I ):

$$\frac{d^2 x(t)}{dt^2} = -\cos\left(\int_0^t cleaf_2(u)du\right)\cdot(cleaf_2(t))^2 - \sin\left(\int_0^t cleaf_2(u)du\right)\cdot\frac{d}{dt}cleaf_2(t) \quad (3.10)$$

where the derivative of $cleaf_2(t)$ with respect to parameter $t$ is as follows:

$$\frac{d}{dt}cleaf_2(t) = -\sqrt{1-(cleaf_2(t))^4} \quad (2m-2)\pi_2 \leq t \leq (2m-1)\pi_2 \quad (3.11)$$

$$\frac{d}{dt}cleaf_2(t) = \sqrt{1-(cleaf_2(t))^4} \quad (2m-1)\pi_2 \leq t \leq 2m\pi_2 \tag{3.12}$$

Notice that the plus or the minus of the derivation: $dcleaf_2(t)/dt$ depend on the range of the parameter $t$[12][13]. The three terms of $d^2x(t)/dt^2$, $x(t)$ and $x(t)^3$ in the Eq. (3.9) are summarized in the table 3. These data are obtained by using the table 2. The value of $d^2x(t)/dt^2 - 3x(t) + 4 x(t)^3$ almost becomes zero as shown in the table 3. The type(Ⅰ) solution is satisfied with the Eq.(3.9) from the numerical data in the table 3.

Table 3 Numerical data in the type (Ⅰ) solution
(All results have been rounded to no more than six significant figures)

| $t$ | $x(t)$ (by the Eq.(3.7)) | $x(t)^3$ (by the Eq.(3.7)) | $d^2x(t)/dt^2$ (by the Eq. (3.10)) | $d^2x(t)/dt^2 - 3x(t) + 4 x(t)^3$ (the Eq.(3.9)) |
|---|---|---|---|---|
| -10.0 | 0.89959 | 0.72801 | -0.21327 | -0.00000 |
| -9.0 | 0.71812 | 0.37033 | 0.67303 | 0.00000 |
| -8.0 | 0.99116 | 0.97372 | -0.92141 | 0.00000 |
| -7.0 | 0.77341 | 0.46264 | 0.46969 | -0.00000 |
| -6.0 | 0.80697 | 0.52550 | 0.31890 | 0.00000 |
| -5.0 | 0.97149 | 0.91689 | -0.75308 | -0.00000 |
| -4.0 | 0.70868 | 0.35593 | 0.70234 | 0.00000 |
| -3.0 | 0.93565 | 0.81913 | -0.46954 | 0.00000 |
| -2.0 | 0.85261 | 0.61981 | 0.07858 | 0.00000 |
| -1.0 | 0.74045 | 0.40597 | 0.59746 | -0.00000 |
| 0.0 | 1.00000 | 1.00000 | -1.00000 | 0.00000 |
| 1.0 | 0.74045 | 0.40597 | 0.59746 | -0.00000 |
| 2.0 | 0.85261 | 0.61981 | 0.07858 | 0.00000 |
| 3.0 | 0.93565 | 0.81913 | -0.46954 | 0.00000 |
| 4.0 | 0.70868 | 0.35593 | 0.70234 | 0.00000 |
| 5.0 | 0.97149 | 0.91689 | -0.75308 | -0.00000 |
| 6.0 | 0.80697 | 0.52550 | 0.31890 | 0.00000 |
| 7.0 | 0.77341 | 0.46264 | 0.46969 | -0.00000 |
| 8.0 | 0.99116 | 0.97372 | -0.92141 | 0.00000 |
| 9.0 | 0.71812 | 0.37033 | 0.67303 | 0.00000 |
| 10.0 | 0.89959 | 0.72801 | -0.21327 | -0.00000 |

The variables are set to $\Phi=0$ and $\omega=1$. The amplitude $A$ is varied. Under these conditions, the wave obtained by the type (I) exact solution is shown in Fig. 3. As the amplitude $A$ is varied, the initial position ( in Eq. (2.4)) also varies. The range of displacement $x(t)$ can be obtained by the following inequality:

$$\frac{1}{\sqrt{2}} A \leq x(t) \leq A \quad (A \geq 0) \tag{3.13}$$

$$A \leq x(t) \leq \frac{1}{\sqrt{2}} A \quad (A < 0) \tag{3.14}$$

The center of the displacement $x(t)$ is obtained as follows:

$$(Center\ of\ displacement) = \frac{1+\frac{1}{\sqrt{2}}}{2} A = \frac{2+\sqrt{2}}{4} A \tag{3.15}$$

The amplitude $A$ is obtained as follows:

$$(Amplitude) = A - \frac{2+\sqrt{2}}{4} A = \frac{2-\sqrt{2}}{4} A \tag{3.16}$$

Next, the variables are set to $\Phi=0$ and $A=1$. The angular frequency $\omega$ is varied. The wave obtained by the type (I) exact solution is shown in Fig. 4. The period of the waves is varied according to the absolute value $\omega$. As the absolute value $\omega$ increases, the period decreases. In contrast, as the absolute value ω decreases, the period increases. In case $\omega=\pm 1$, the period becomes the constant $\pi_2$. In case $\omega=\pm 2$, the period becomes the constant $\pi_2/2$. In case $\omega=\pm 3$, the period becomes the constant $\pi_2/3$. By using the parameter $\omega$, the period $T$ is obtained as follows:

$$T = \frac{\pi_2}{|\omega|} \tag{3.17}$$

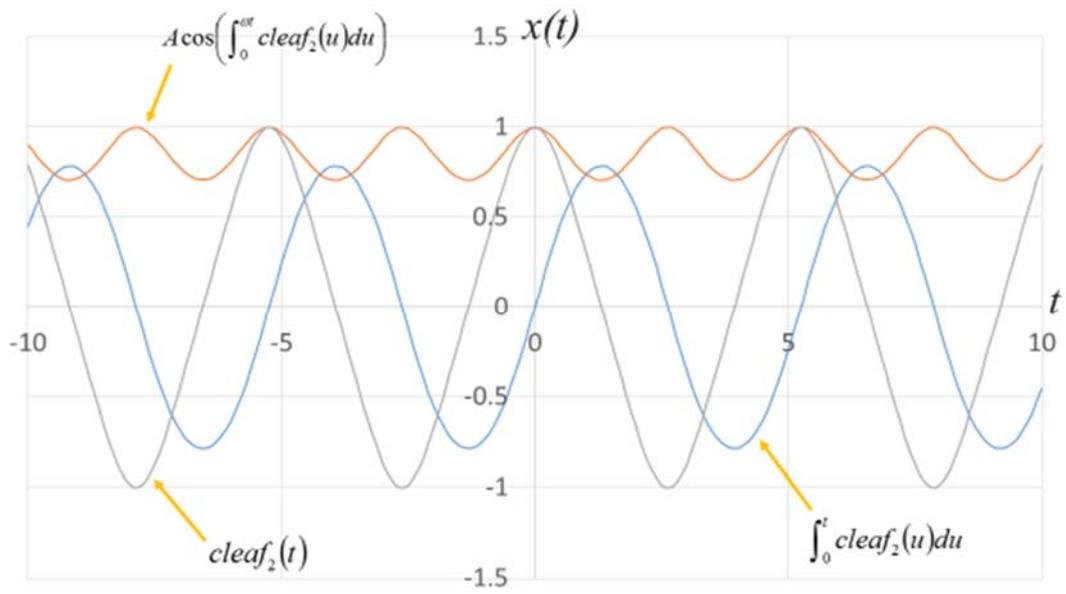

Fig. 2 Waves of the exact solution obtained by Type(I); the leaf function: $cleaf_2(t)$ and the function: $\int_0^t cleaf_2(u)du$

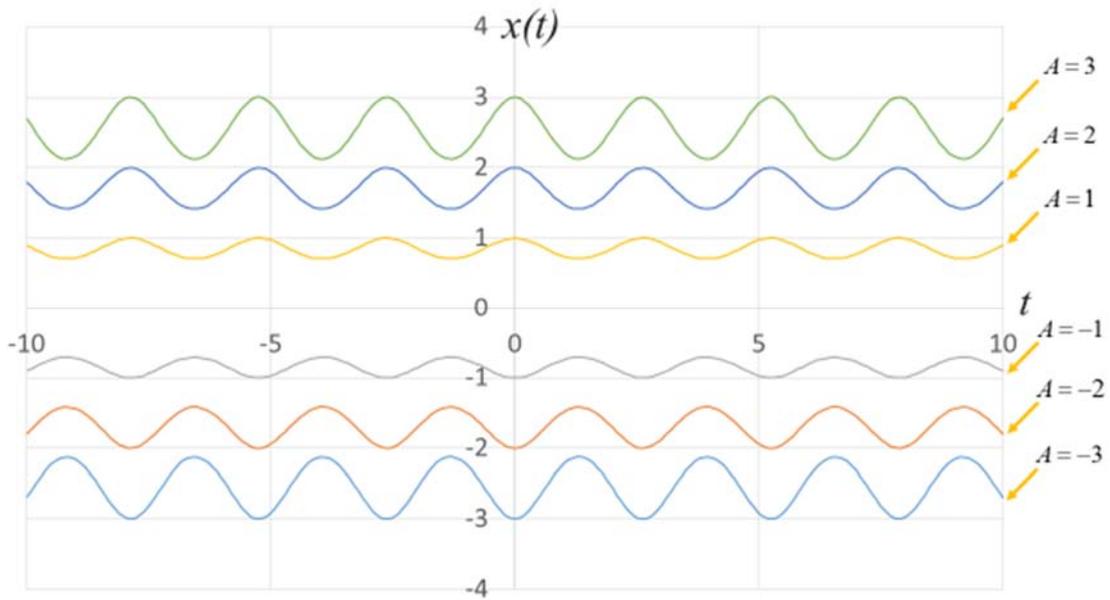

Fig. 3 Wave obtained by Type(I) exact solution with respect to the variation in the amplitude $A$

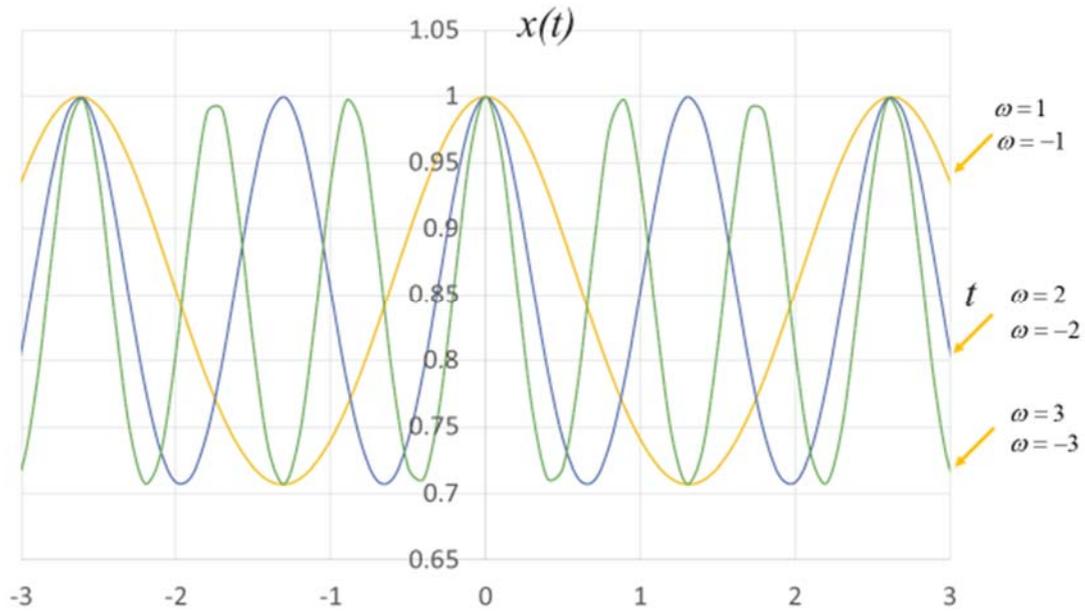

Fig. 4 Wave obtained by Type(I) exact solution with respect to the variation in the angular frequency $\omega$

### 3.2 Numerical Results of Exact Solution of Type(Ⅱ)

In the case $A=1$, $\omega=1$, $\Phi=0$ in Eq. (2.6), the exact solution of type (Ⅱ), the second derivative and the ordinary difference equation are as follows:

$$x(t) = \sin\left(\int_0^t cleaf_2(u)du\right) \tag{3.18}$$

$$\frac{d^2x(t)}{dt^2} = -\sin\left(\int_0^t cleaf_2(u)du\right)\cdot(cleaf_2(t))^2 + \cos\left(\int_0^t cleaf_2(u)du\right)\cdot\frac{d}{dt}cleaf_2(t) \tag{3.19}$$

$$\frac{d^2x(t)}{dt^2} + 3x(t) - 4x(t)^3 = 0 \tag{3.20}$$

The data of the terms: $d^2x(t)/dt^2$, $x(t)$ and $x(t)^3$ in the Eq. (3.20) are summarized in the table 4. The type(Ⅱ) solution is satisfied with the Eq.(3.20), as shown in the table 4.

Table 4 Numerical data in the type (II) solution
(All results have been rounded to no more than six significant figures)

| t | x(t) (by the Eq.(3.18)) | x(t)³ (by the Eq.(3.18)) | d²x(t)/dt² (by the Eq. (3.19)) | d²x(t)/dt²+3x(t)- 4 x(t)³ (the Eq.(3.20)) |
|---|---|---|---|---|
| -10.0 | 0.43672 | 0.08329 | -0.97699 | 0.00000 |
| -9.0 | 0.69591 | 0.33703 | -0.73961 | 0.00000 |
| -8.0 | 0.13263 | 0.00233 | -0.38858 | -0.00000 |
| -7.0 | -0.63389 | -0.25471 | 0.88283 | 0.00000 |
| -6.0 | -0.59059 | -0.20599 | 0.94778 | -0.00000 |
| -5.0 | 0.23707 | 0.01332 | -0.65791 | 0.00000 |
| -4.0 | 0.70552 | 0.35118 | -0.71184 | 0.00000 |
| -3.0 | 0.35290 | 0.04395 | -0.88290 | -0.00000 |
| -2.0 | -0.52253 | -0.14267 | 0.99690 | -0.00000 |
| -1.0 | -0.67210 | -0.30360 | 0.80189 | 0.00000 |
| 0.0 | 0.00000 | 0.00000 | 0.00000 | 0.00000 |
| 1.0 | 0.67210 | 0.30360 | -0.80189 | -0.00000 |
| 2.0 | 0.52253 | 0.14267 | -0.99690 | 0.00000 |
| 3.0 | -0.35290 | -0.04395 | 0.88290 | 0.00000 |
| 4.0 | -0.70552 | -0.35118 | 0.71184 | -0.00000 |
| 5.0 | -0.23707 | -0.01332 | 0.65791 | -0.00000 |
| 6.0 | 0.59059 | 0.20599 | -0.94778 | 0.00000 |
| 7.0 | 0.63389 | 0.25471 | -0.88283 | -0.00000 |
| 8.0 | -0.13263 | -0.00233 | 0.38858 | 0.00000 |
| 9.0 | -0.69591 | -0.33703 | 0.73961 | -0.00000 |
| 10.0 | -0.43672 | -0.08329 | 0.97699 | -0.00000 |

In the case, the wave of the exact solution of type (II) is compared to the wave of the function: $cleaf_2(t)$ and the wave of the function: $\int_0^t cleaf_2(u)du$. These waves are shown in Fig. 5. The horizontal and vertical axes represent time and displacement *x(t)*, respectively. The maximum value *x(t)* of the Type(II) exact solution is obtained as follows:

$$\sin\left(\int_0^{\frac{\pi_2}{2}(4m+1)} cleaf_2(u)du\right) = \sin\left(\frac{\pi}{4}\right) = \frac{1}{\sqrt{2}} \quad (m: integer) \quad (3.21)$$

The minimum value $x(t)$ of the Type(II) exact solution is obtained as follows:

$$\sin\left(\int_0^{\frac{\pi_2}{2}(4m-1)} cleaf_2(u)du\right) = \sin\left(-\frac{\pi}{4}\right) = -\frac{1}{\sqrt{2}} \quad (m: integer) \tag{3.22}$$

The range of the variable $x(t)$ becomes $-\frac{1}{\sqrt{2}} \leq x(t) \leq \frac{1}{\sqrt{2}}$. The center of the displacement and amplitude become $0.0$ and $\frac{1}{\sqrt{2}}$, respectively. The period of the type (II) exact solution becomes $2\pi_2$.

Next, the variables are set to $\Phi=0$ and $\omega=1$. The amplitude $A$ is varied. Under these conditions, the wave obtained by the type (II) exact solution is shown in Figs. 6 and 7. The center of the displacement $x(t)$ becomes 0.0. The range of displacement $x(t)$ can be obtained by the following inequality:

$$-\frac{1}{\sqrt{2}}|A| \leq x(t) \leq \frac{1}{\sqrt{2}}|A| \tag{3.23}$$

The amplitude is obtained as follows:

$$(Amplitude) = \frac{1}{\sqrt{2}}|A| \tag{3.24}$$

Next, the variables are set to $\Phi=0$ and $A=1$. The angular frequency $\omega$ is varied. The wave obtained by the type (II) exact solution is shown in Figs. 8 and 9. In case $\omega=\pm1$, the period becomes the constant $2\pi_2$. In case $\omega=\pm2$, the period becomes the constant $2\pi_2/2$. In case $\omega=\pm3$, the period becomes the constant $2\pi_2/3$. By using the parameter $\omega$, the period $T$ is obtained as follows:

$$T = \frac{2\pi_2}{|\omega|} \tag{3.25}$$

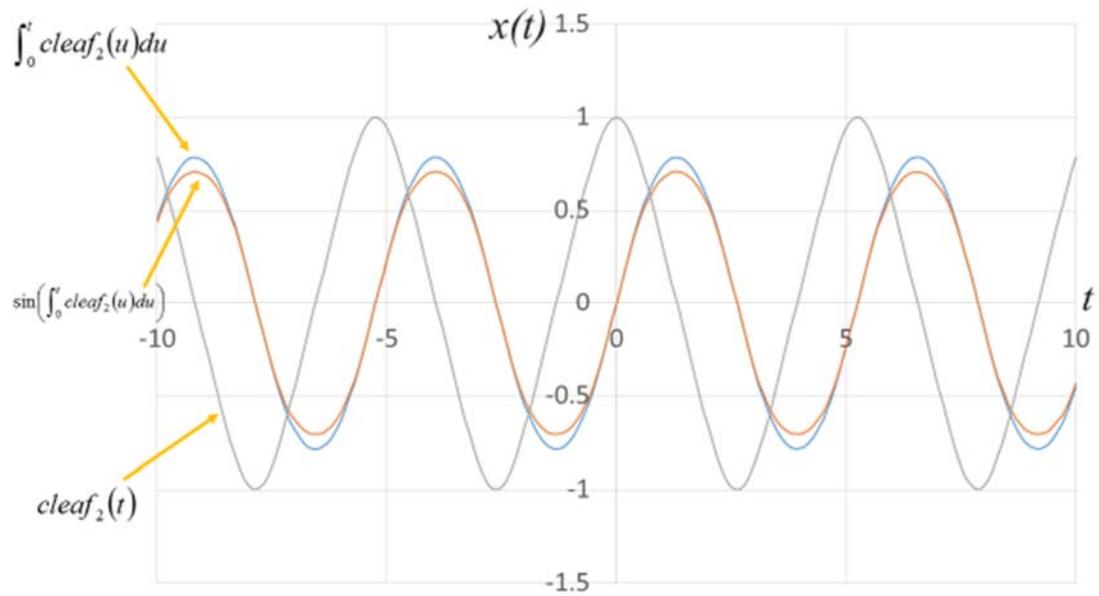

Fig. 5 Waves of the exact solution obtained by Type(II); the leaf function: $cleaf_2(t)$ and the function: $\int_0^t cleaf_2(u)du$

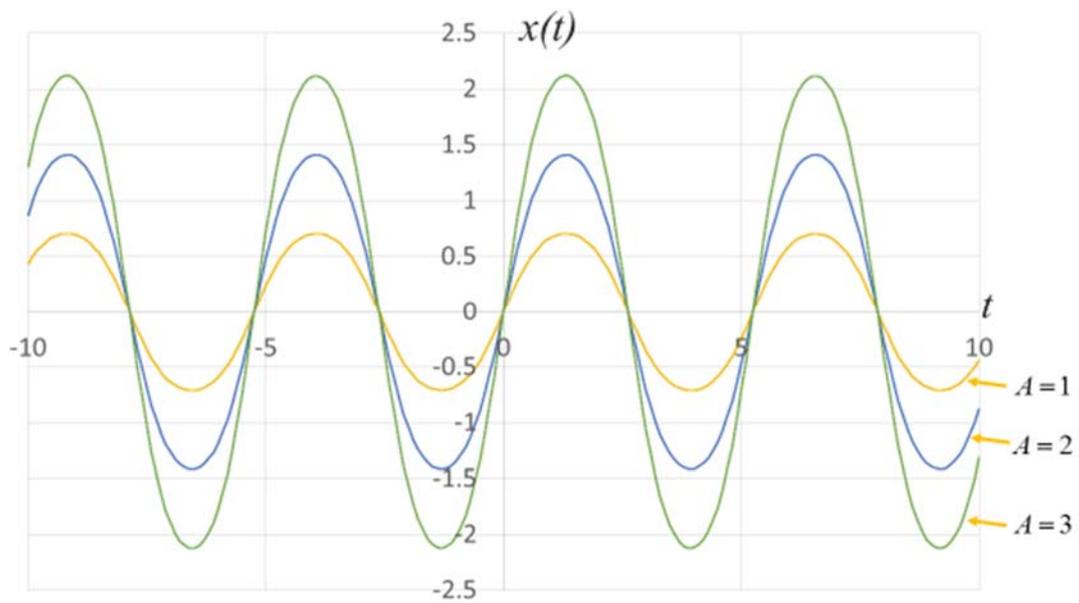

Fig. 6 Wave obtained by Type(II) exact solution with respect to the variation in the amplitude A ($A=1, 2, 3$)

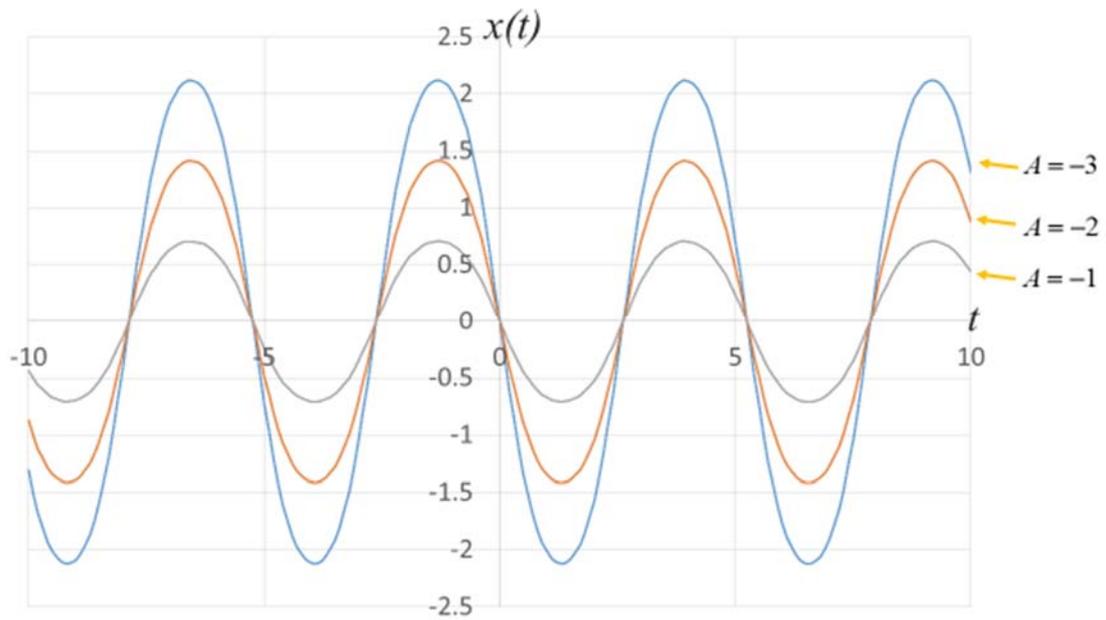

Fig. 7 Wave obtained by Type(II) exact solution with respect to the variation in the amplitude A (*A=-1, -2, -3*)

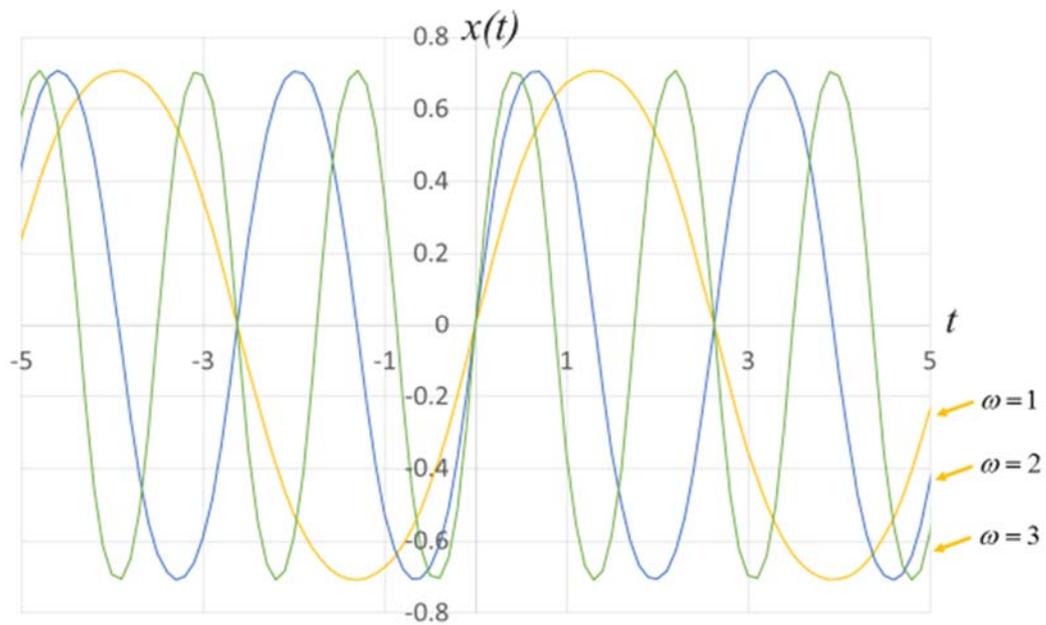

Fig. 8 Wave obtained by Type (II) exact solution with respect to the variation in the angular frequency $\omega$ (*$\omega$=1, 2, 3*)

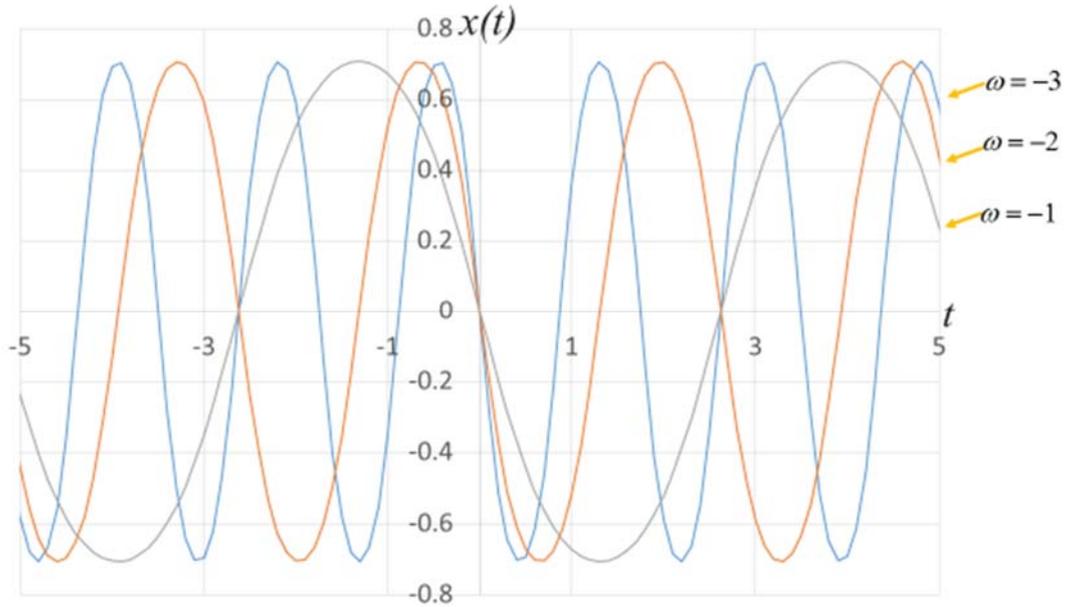

Fig. 9 Wave obtained by Type (II) exact solution with respect to the variation in the angular frequency ω (ω=-1, -2, -3 )

### 3.3 Numerical Results of Exact Solution of Type(Ⅲ)

In the case $A=1$, $\omega=1$, $\Phi=0$ in Eq. (2.10), the exact solution of type (Ⅲ), the second derivative and the ordinary difference equation are as follows:

$$x(t) = \cos\left(\int_0^t sleaf_2(u)du\right) + \sin\left(\int_0^t sleaf_2(u)du\right) \qquad (3.26)$$

$$\frac{d^2x(t)}{dt^2} = 4\left\{\cos\left(\int_0^t sleaf_2(u)du\right)\right\}^3 + 4\left\{\sin\left(\int_0^t sleaf_2(u)du\right)\right\}^3 - 3\cos\left(\int_0^t sleaf_2(u)du\right) - 3\sin\left(\int_0^t sleaf_2(u)du\right) \qquad (3.27)$$

$$\frac{d^2x(t)}{dt^2} - 3x(t) + 2x(t)^3 = 0 \qquad (3.28)$$

The data of the terms: $d^2x(t)/dt^2$, $x(t)$ and $x(t)^3$ in the Eq. (3.28) are summarized in the table 5. The type(Ⅲ) solution is satisfied with the Eq.(3.28) from the numerical data in the table 5.

$$x(t) = \cos\left(\int_0^t sleaf_2(u)du\right) + \sin\left(\int_0^t sleaf_2(u)du\right) \tag{3.29}$$

Table 5 Numerical data in the type (III) solution
(All results have been rounded to no more than six significant figures)

| $t$ | $\cos\left(\int_0^t sleaf_2(u)du\right)$ | $\sin\left(\int_0^t sleaf_2(u)du\right)$ | $x(t)$ (by the Eq.(3.26)) | $x(t)^3$ (by the Eq.(3.26)) | $d^2x(t)/dt^2$ (by the Eq. (3.27)) | $d^2x(t)/dt^2 - 3x(t) + 2 x(t)^3$ (the Eq.(3.28)) |
|---|---|---|---|---|---|---|
| -10.0 | 0.99293 | 0.11867 | 1.11161 | 1.37359 | 0.58764 | 0.00000 |
| -9.0 | 0.57289 | 0.81962 | 1.39252 | 2.70027 | -1.22297 | 0.00000 |
| -8.0 | 0.00895 | 0.99995 | 1.00891 | 1.02698 | 0.97277 | 0.00000 |
| -7.0 | 0.36000 | 0.93294 | 1.29295 | 2.16148 | -0.44410 | 0.00000 |
| -6.0 | 0.96032 | 0.27887 | 1.23920 | 1.90294 | -0.08828 | 0.00000 |
| -5.0 | 0.99955 | 0.02978 | 1.02934 | 1.09064 | 0.90674 | 0.00000 |
| -4.0 | 0.75273 | 0.65831 | 1.41105 | 2.80953 | -1.38589 | 0.00000 |
| -3.0 | 0.07131 | 0.99745 | 1.06876 | 1.22082 | 0.76466 | 0.00000 |
| -2.0 | 0.19132 | 0.98152 | 1.17285 | 1.61336 | 0.29183 | 0.00000 |
| -1.0 | 0.88508 | 0.46542 | 1.35051 | 2.46318 | -0.87483 | 0.00000 |
| 0.0 | 1.00000 | 0.00000 | 1.00000 | 1.00000 | 1.00000 | 0.00000 |
| 1.0 | 0.88508 | 0.46542 | 1.35051 | 2.46318 | -0.87483 | 0.00000 |
| 2.0 | 0.19132 | 0.98152 | 1.17285 | 1.61336 | 0.29183 | 0.00000 |
| 3.0 | 0.07131 | 0.99745 | 1.06876 | 1.22082 | 0.76466 | 0.00000 |
| 4.0 | 0.75273 | 0.65831 | 1.41105 | 2.80953 | -1.38589 | 0.00000 |
| 5.0 | 0.99955 | 0.02978 | 1.02934 | 1.09064 | 0.90674 | 0.00000 |
| 6.0 | 0.96032 | 0.27887 | 1.23920 | 1.90294 | -0.08828 | 0.00000 |
| 7.0 | 0.36000 | 0.93294 | 1.29295 | 2.16148 | -0.44410 | 0.00000 |
| 8.0 | 0.00895 | 0.99995 | 1.00891 | 1.02698 | 0.97277 | 0.00000 |
| 9.0 | 0.57289 | 0.81962 | 1.39252 | 2.70027 | -1.22297 | 0.00000 |
| 10.0 | 0.99293 | 0.11867 | 1.11161 | 1.37359 | 0.58764 | 0.00000 |

In the case $A=1$, $\omega=1$, $\Phi=0$ in Eq. (2.10), the wave of the exact solution of type (III) is compared to the wave of the function: $sleaf_2(t)$ and the wave of the function: $\int_0^t sleaf_2(u)du$. These waves are shown in Fig. 10. The horizontal and vertical axes represent time and displacement $x(t)$, respectively. In the function: $sleaf_2(t)$, the amplitude and the period become 1.0 and $2\pi_2(\cong 5.244)$[12][13], respectively. The center of the displacement becomes $x(t)=0$. In the function: $\int_0^t sleaf_2(u)du$, the minimum value of the displacement becomes 0.0. The maximum value of the displacement is given as follows:

$$\int_0^{\pi_2} sleaf_2(u)du = \frac{\pi}{2}(\cong 1.5708) \tag{3.30}$$

The center of the displacement is obtained as follows:

$$(Center\ of\ displacement) = \frac{0+\frac{\pi}{2}}{2} = \frac{\pi}{4}(\cong 0.7854\cdots) \tag{3.31}$$

The amplitude of the function: $\int_0^t sleaf_2(u)du$ becomes $\pi/4$.

Next, the wave obtained by the type(Ⅲ) exact solution is discussed. By using the addition theorem, the type(Ⅲ) exact solution can be transformed as follows:

$$x(t) = \sqrt{2}\sin\left(\int_0^t sleaf_2(u)du + \frac{\pi}{4}\right) \tag{3.32}$$

The minimum of the variable $x(t)$ is obtained as follows:

$$x(2m\pi_2) = \sqrt{2}\sin\left(\int_0^{2m\pi_2} sleaf_2(u)du + \frac{\pi}{4}\right) = \sqrt{2}\sin\left(0 + \frac{\pi}{4}\right) = 1.0 \quad (m:\ integer) \tag{3.33}$$

or

$$x((2m-1)\pi_2) = \sqrt{2}\sin\left(\int_0^{(2m-1)\pi_2} sleaf_2(u)du + \frac{\pi}{4}\right) = \sqrt{2}\sin\left(\frac{\pi}{2} + \frac{\pi}{4}\right) = 1.0 \quad (m:integer) \tag{3.34}$$

In contrast, the maximum value is obtained as follows:

$$\sqrt{2}\sin\left(\int_0^{(2m-1)\frac{\pi_2}{2}} sleaf_2(u)du + \frac{\pi}{4}\right) = \sqrt{2}\sin\left(\frac{\pi}{4} + \frac{\pi}{4}\right) = \sqrt{2} \quad (m:\ integer) \tag{3.35}$$

The range of the variable $x(t)$ becomes $1.0 \leq x(t) \leq \sqrt{2}$. The center of the displacement and amplitude becomes $x(t) = \frac{1+\sqrt{2}}{2}$ and $\frac{1+\sqrt{2}}{2}$, respectively.

The variables are set to $\Phi=0$ and $\omega=1$. The amplitude $A$ is varied. Under these conditions, the wave obtained by the type (III) exact solution is shown in Fig. 11. As the amplitude A is varied, the initial position (in Eq. (2.12)) also varies. The range of the displacement $x(t)$ can be obtained by the following inequality:

$$A \le x(t) \le \sqrt{2}A \quad (A \ge 0) \tag{3.36}$$

$$\sqrt{2}A \le x(t) \le A \quad (A < 0) \tag{3.37}$$

The center of the displacement $x(t)$ is obtained as follows:

$$(Center\ of\ displacement) = \frac{1+\sqrt{2}}{2}A \tag{3.38}$$

The amplitude is obtained as follows:

$$(Amplitude) = \frac{1+\sqrt{2}}{2}|A| \tag{3.39}$$

Next, the variables are set to $\Phi=0$ and $A=1$. The angular frequency $\omega$ is varied. The wave obtained by the type (III) exact solution is shown in Fig. 12. In case $\omega=\pm 1$, the period becomes the constant $\pi_2$. In case $\omega=\pm 2$, the period becomes the constant $\pi_2/2$. In case $\omega=\pm 3$, the period becomes the constant $\pi_2/3$. By using the parameter $\omega$, the period $T$ is obtained as follows:

$$T = \frac{\pi_2}{|\omega|} \tag{3.40}$$

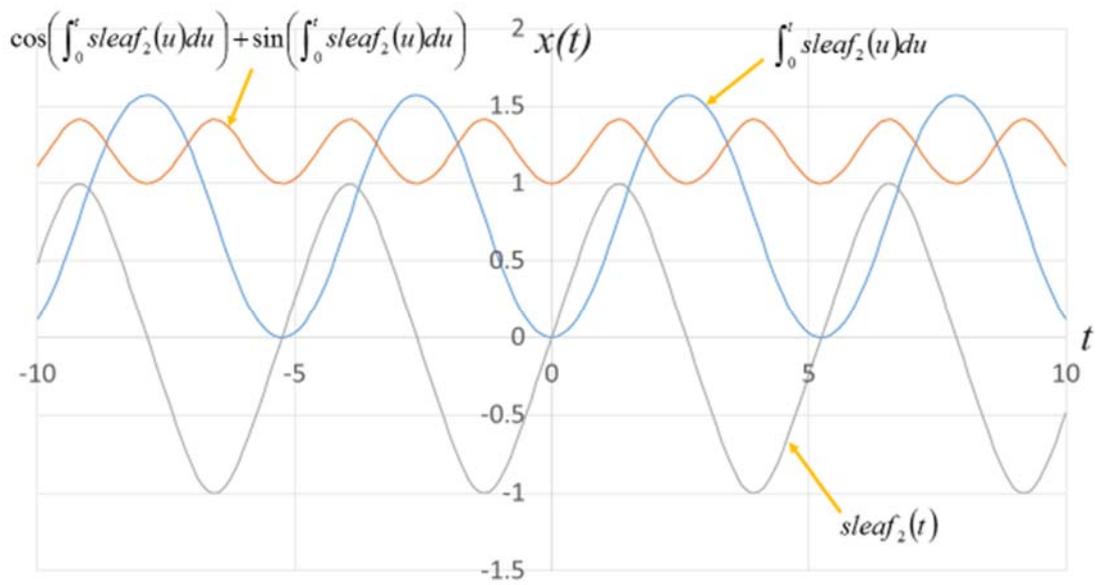

Fig. 10 Waves of the exact solution obtained by Type(III); the leaf function: $sleaf_2(t)$ and the function: $\int_0^t sleaf_2(u)du$

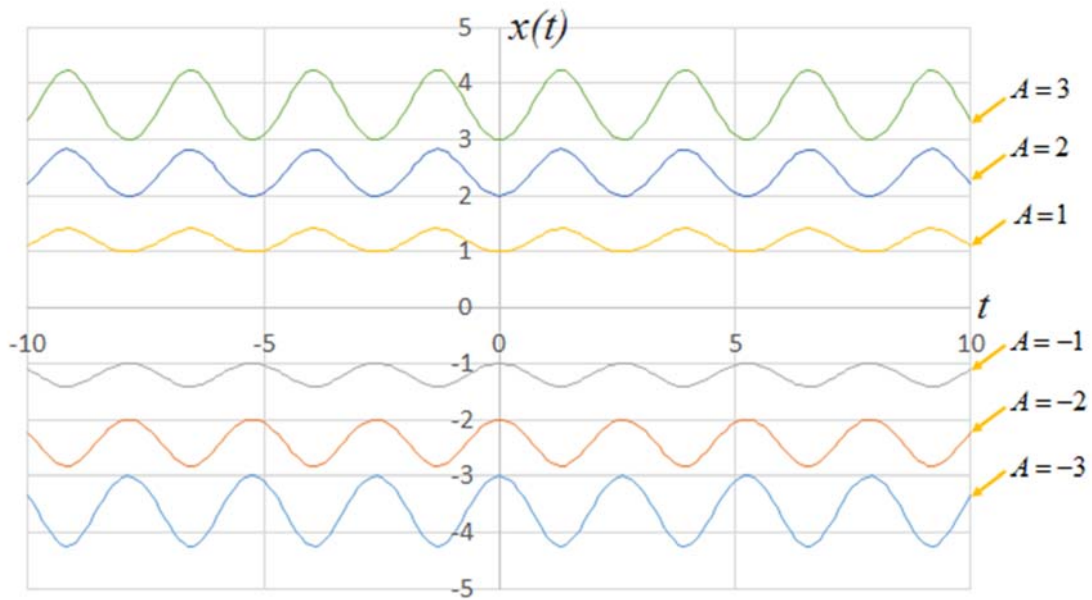

Fig. 11 Wave obtained by Type(III) exact solution with respect to the variation in the amplitude $A$

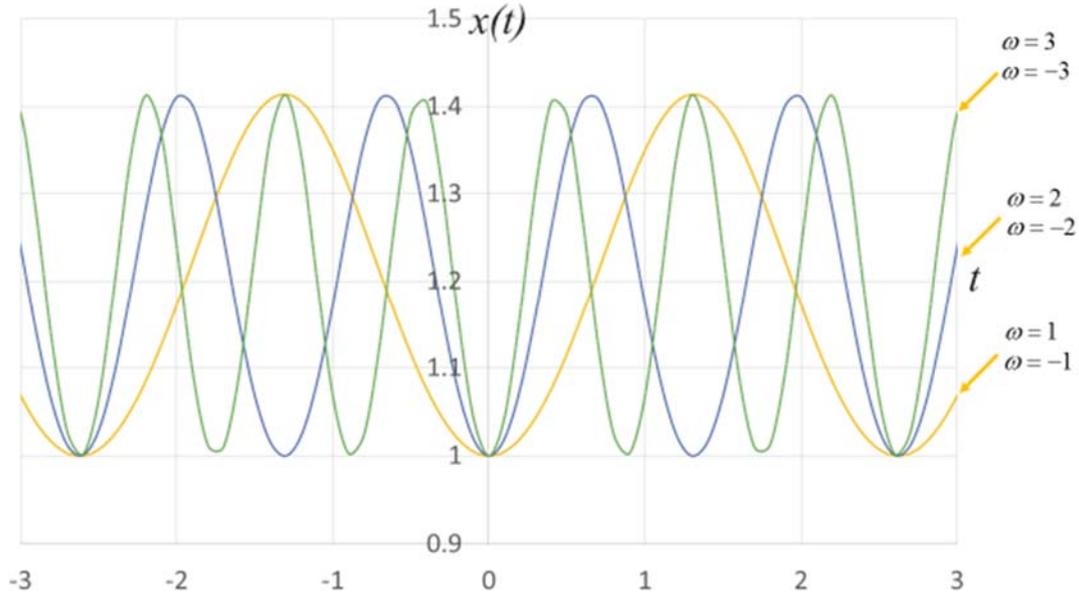

Fig. 12 Wave obtained by Type (III) exact solution with respect to the variation in the angular frequency $\omega$

### 3.4 Numerical Results of Exact Solution of Type(IV)

In the case $A=1,\ \omega=1,\ \Phi=0$ in Eq. (2.14), the wave of the exact solution of type (IV), the wave of the function: $sleaf_2(t)$, and the wave of the function: $\int_0^t sleaf_2(u)du$ are shown in Fig. 13. The horizontal and vertical axes represent time and displacement $x(t)$, respectively. The numerical data of the type (IV) solution can be obtained by using the table 5. By using the addition theorem, the type(IV) exact solution can be transformed as follows:

$$x(t) = -\sqrt{2}\sin\left(\int_0^t sleaf_2(u)du - \frac{\pi}{4}\right) \qquad (3.41)$$

The maximum value of the displacement is given as follows:

$$x(2m\pi_2) = -\sqrt{2}\sin\left(\int_0^{2m\pi_2} sleaf_2(u)du - \frac{\pi}{4}\right) = -\sqrt{2}\sin\left(0 - \frac{\pi}{4}\right) = 1.0\ (m\text{:}integer) \qquad (3.42)$$

In contrast, the minimum value of the displacement is given as follows:

$$x((2m-1)\pi_2) = -\sqrt{2}\sin\left(\int_0^{(2m-1)\pi_2} sleaf_2(u)du - \frac{\pi}{4}\right) = -\sqrt{2}\sin\left(\frac{\pi}{2} - \frac{\pi}{4}\right) = -1.0 \;(m: integer) \quad (3.43)$$

The range of the variable $x(t)$ becomes $-1.0 \leqq x(t) \leqq 1.0$. The center of the displacement and amplitude becomes $x(t)=0$ and $1.0$, respectively.

Next, the variables are set to $\Phi=0$ and $\omega=1$. The amplitude $A$ is varied. Under these conditions, the wave obtained by the type (IV) exact solution is shown in Figs. 14 and 15. The range of the displacement $x(t)$ can be obtained by the following inequality:

$$-|A| \leq x(t) \leq |A| \qquad (3.44)$$

The center of the displacement and amplitude become $x(t)=0$ and $|A|$, respectively. Next, the variables are set to $\Phi=0$ and $A=1$. The angular frequency $\omega$ is varied. The wave obtained by the type (IV) exact solution is shown in Fig. 16. The period of the waves is varied according to the absolute value $\omega$. As the absolute value $\omega$ increases, the period decreases. In contrast, as the absolute value $\omega$ decreases, the period increases. In case $\omega=\pm 1$, the period becomes the constant $2\pi_2$. In case $\omega=\pm 2$, the period becomes the constant $2\pi_2/2$. In case $\omega=\pm 3$, the period becomes the constant $2\pi_2/3$. By using the parameter $\omega$, the period $T$ is obtained as follows:

$$T = \frac{2\pi_2}{|\omega|} \qquad (3.45)$$

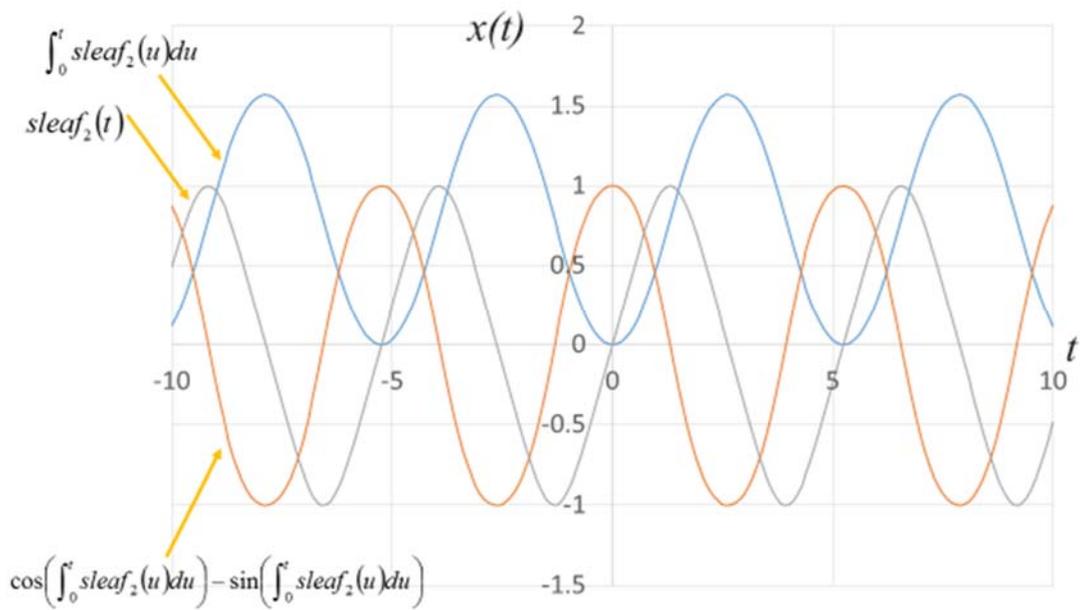

Fig. 13 Waves of the exact solution obtained by Type(IV); the leaf function: $sleaf_2(t)$ and the function: $\int_0^t sleaf_2(u)du$

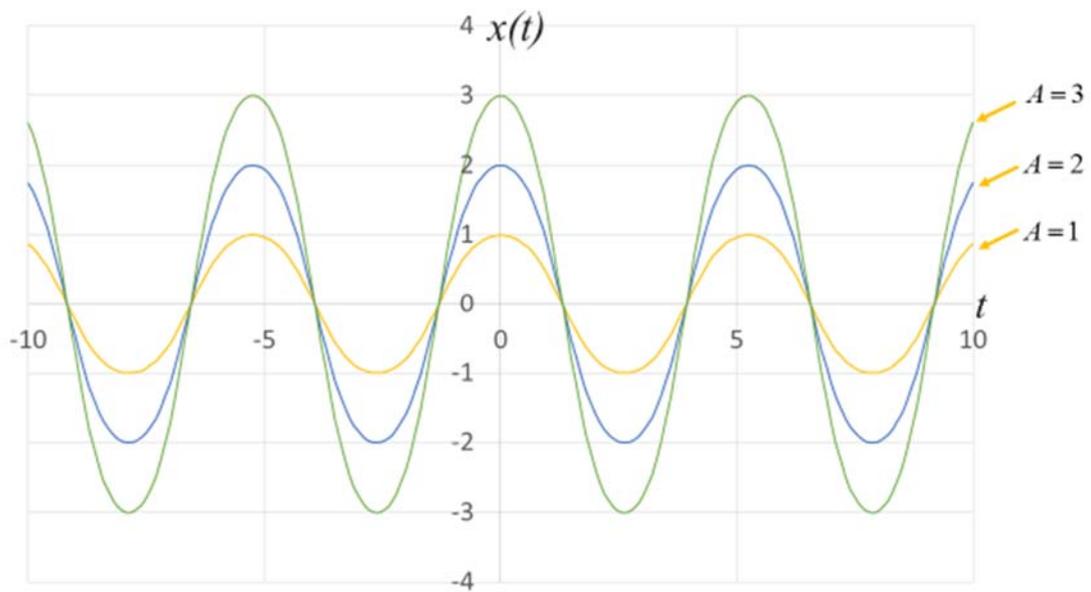

Fig. 14 Wave obtained by Type(IV) exact solution with respect to the variation in the amplitude $A$ ($A=1, 2, 3$)

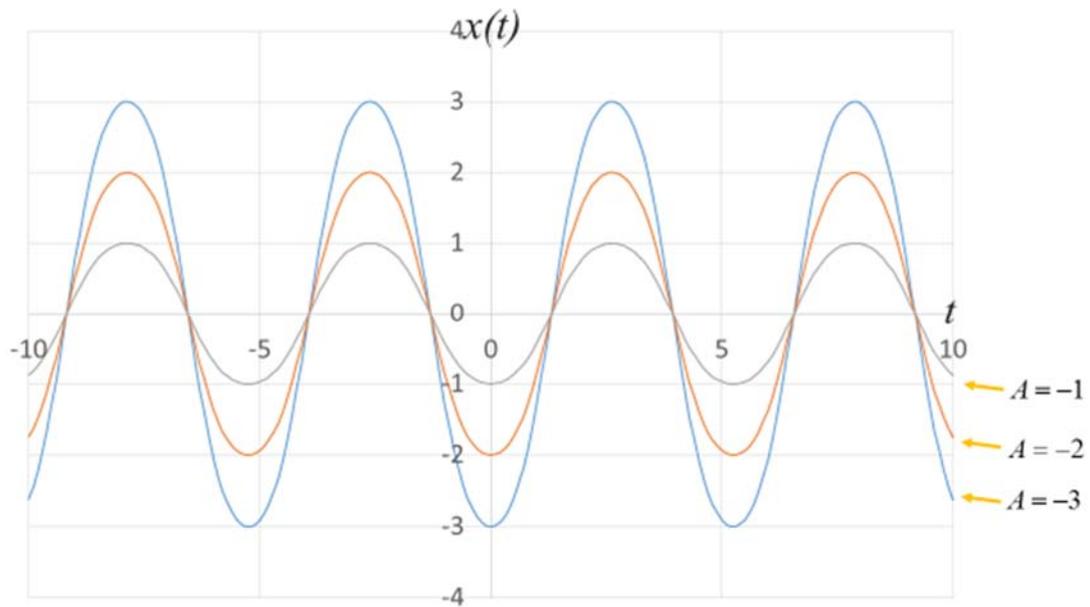

Fig. 15 Wave obtained by Type(IV) exact solution with respect to the variation in the amplitude $A$ ($A=-1, -2, -3$)

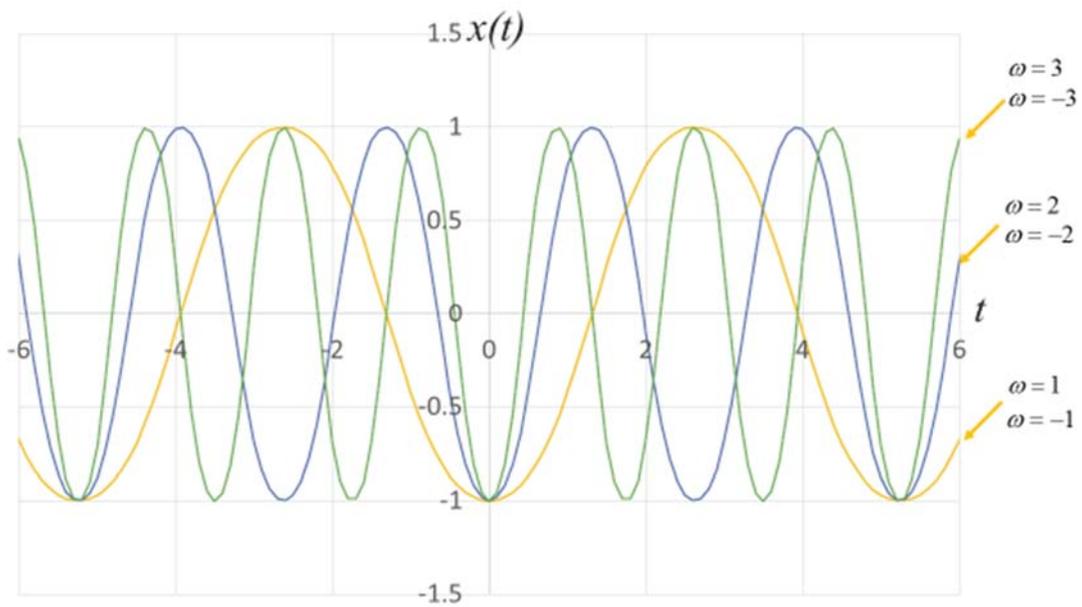

Fig. 16 Wave obtained by Type (IV) exact solution with respect to the variation in the angular frequency $\omega$

### 3.5. Numerical Results of Exact Solution of Type( V)

In the case $A=1$, $\omega=1$, $\Phi=0$ in Eq. (2.18), the wave of the exact solution of type (V), the wave of the leaf function: $sleaf_2(t)$ and $cleaf_2(t)$, and the wave of the function: $\int_0^t sleaf_2(u)du$ and $\int_0^t cleaf_2(u)du$ are shown are shown in Fig. 17. In order to show the wave of the type(V) exact solution, Fig. 18 shows an enlarged view of Fig. 17. The horizontal and vertical axes represent time $t$ and displacement $x(t)$, respectively. The numerical data of the Type (V) solution can be obtained by using the table 3 and the table 5. To discuss the range of the type(V) exact solution $x(t)$, the following equation is transformed:

$$x_1(t) = \cos\left(\int_0^t sleaf_2(u)du\right) + \sin\left(\int_0^t sleaf_2(u)du\right) \tag{3.46}$$

The following equation is obtained by squaring both sides of the above equation.

$$\begin{aligned}\{x_1(t)\}^2 &= \left\{\cos\left(\int_0^t sleaf_2(u)du\right)\right\}^2 + \left\{\sin\left(\int_0^t sleaf_2(u)du\right)\right\}^2 \\ &+ 2\sin\left(\int_0^t sleaf_2(u)du\right)\cos\left(\int_0^t sleaf_2(u)du\right) \\ &= 1 + \sin\left(2\int_0^t sleaf_2(u)du\right) = 1 + (sleaf_2(t))^2\end{aligned} \tag{3.47}$$

Eq. (Ⅲ.4) in the appendix Ⅲ is applied to the above equation. As shown in Fig. 10, the inequality: $x_1(t)>0$ is obvious. The following equation is obtained.

$$x_1(t) = \sqrt{1 + (sleaf_2(t))^2} \tag{3.48}$$

Next, the following equation is transformed:

$$x_2(t) = \sqrt{2}\cos\left(\int_0^t cleaf_2(u)du\right) \tag{3.49}$$

The following equation is obtained by squaring both sides of the above equation.

$$\{x_2(t)\}^2 = 2\left\{\cos\left(\int_0^t cleaf_2(u)du\right)\right\}^2 = 2\frac{1+\cos\left(2\int_0^t cleaf_2(u)du\right)}{2}$$
$$= 1+\cos\left(2\int_0^t cleaf_2(u)du\right) = 1+(cleaf_2(t))^2 \qquad (3.50)$$

As shown in Fig. 2, the inequality: $x_2(t)>0$ is obvious. The following equation is obtained.

$$x_2(t) = \sqrt{1+(cleaf_2(t))^2} \qquad (3.51)$$

Therefore, the type(V) exact solution is obtained as follows:

$$x(t) = x_1(t) + x_2(t) = \sqrt{1+(sleaf_2(t))^2} + \sqrt{1+(cleaf_2(t))^2} \qquad (3.52)$$

The sign of the first order differential with respect to the leaf function $sleaf_2(t)$ and $cleaf_2(t)$ depends on the domain $t$ of the variable $x(t)$ (See Ref. [12][13] or the Eqs. (3.11)-(3.12)). The first order differential of equation (3.52) is discussed by division into four domains (1)–(4).

Domain (1): $2m\pi_2 \leq t < \frac{\pi_2}{2} + 2m\pi_2$ ($m$ :integer)

$$\frac{dx(t)}{dt} = \frac{d}{dt}\left(1+(sleaf_2(t))^2\right)^{\frac{1}{2}} + \frac{d}{dt}\left(1+(cleaf_2(t))^2\right)^{\frac{1}{2}}$$
$$= \frac{1}{2}\left(1+(sleaf_2(t))^2\right)^{\frac{1}{2}-1} \cdot 2(sleaf_2(t)) \cdot \sqrt{1-(sleaf_2(t))^4}$$
$$+ \frac{1}{2}\left(1+(cleaf_2(t))^2\right)^{\frac{1}{2}-1} \cdot 2(cleaf_2(t)) \cdot \left(-\sqrt{1-(cleaf_2(t))^4}\right) \qquad (3.53)$$
$$= \frac{(sleaf_2(t))\sqrt{(1+(sleaf_2(t))^2)(1-(sleaf_2(t))^2)}}{\sqrt{1+(sleaf_2(t))^2}} - \frac{(cleaf_2(t))\sqrt{(1+(cleaf_2(t))^2)(1-(cleaf_2(t))^2)}}{\sqrt{1+(cleaf_2(t))^2}}$$
$$= sleaf_2(t)\sqrt{1-(sleaf_2(t))^2} - cleaf_2(t)\sqrt{1-(cleaf_2(t))^2}$$

Domain (2): $\frac{\pi_2}{2} + 2m\pi_2 \leq t < \pi_2 + 2m\pi_2$ ($m$ :integer)

$$\frac{dx(t)}{dt} = \frac{d}{dt}\left(1+(sleaf\,_2(t))^2\right)^{\frac{1}{2}} + \frac{d}{dt}\left(1+(cleaf\,_2(t))^2\right)^{\frac{1}{2}}$$
$$= \frac{1}{2}\left(1+(sleaf\,_2(t))^2\right)^{\frac{1}{2}-1} \cdot 2(sleaf\,_2(t)) \cdot \left(-\sqrt{1-(sleaf\,_2(t))^4}\right)$$
$$+ \frac{1}{2}\left(1+(cleaf\,_2(t))^2\right)^{\frac{1}{2}-1} \cdot 2(cleaf\,_2(t)) \cdot \left(-\sqrt{1-(cleaf\,_2(t))^4}\right)$$
$$= -sleaf\,_2(t)\sqrt{1-(sleaf\,_2(t))^2} - cleaf\,_2(t)\sqrt{1-(cleaf\,_2(t))^2}$$
(3.54)

Domain (3): $\pi_2 + 2m\pi_2 \le t < \frac{3\pi_2}{2} + 2m\pi_2$ (m :integer)

$$\frac{dx(t)}{dt} = \frac{d}{dt}\left(1+(sleaf\,_2(t))^2\right)^{\frac{1}{2}} + \frac{d}{dt}\left(1+(cleaf\,_2(t))^2\right)^{\frac{1}{2}}$$
$$= \frac{1}{2}\left(1+(sleaf\,_2(t))^2\right)^{\frac{1}{2}-1} \cdot 2(sleaf\,_2(t)) \cdot \left(-\sqrt{1-(sleaf\,_2(t))^4}\right)$$
$$+ \frac{1}{2}\left(1+(cleaf\,_2(t))^2\right)^{\frac{1}{2}-1} \cdot 2(cleaf\,_2(t)) \cdot \left(\sqrt{1-(cleaf\,_2(t))^4}\right)$$
$$= -sleaf\,_2(t)\sqrt{1-(sleaf\,_2(t))^2} + cleaf\,_2(t)\sqrt{1-(cleaf\,_2(t))^2}$$
(3.55)

Domain (4): $\frac{3\pi_2}{2} + 2m\pi_2 \le t < 2\pi_2 + 2m\pi_2$ (m :integer)

$$\frac{dx(t)}{dt} = \frac{d}{dt}\left(1+(sleaf\,_2(t))^2\right)^{\frac{1}{2}} + \frac{d}{dt}\left(1+(cleaf\,_2(t))^2\right)^{\frac{1}{2}}$$
$$= \frac{1}{2}\left(1+(sleaf\,_2(t))^2\right)^{\frac{1}{2}-1} \cdot 2(sleaf\,_2(t)) \cdot \left(\sqrt{1-(sleaf\,_2(t))^4}\right)$$
$$+ \frac{1}{2}\left(1+(cleaf\,_2(t))^2\right)^{\frac{1}{2}-1} \cdot 2(cleaf\,_2(t)) \cdot \left(\sqrt{1-(cleaf\,_2(t))^4}\right)$$
$$= sleaf\,_2(t)\sqrt{1-(sleaf\,_2(t))^2} + cleaf\,_2(t)\sqrt{1-(cleaf\,_2(t))^2}$$
(3.56)

The extreme value of the type(V) solution is obtained by the equation: *dx(t)/dt=0* with respect to Eqs. (3.53)–(3.56). By using the equation: *dx(t)/dt=0*, the same following equation is derived by Eqs. (3.53)–(3.56).

$$(sleaf\,_2(t) - cleaf\,_2(t)) \cdot (sleaf\,_2(t) + cleaf\,_2(t)) \cdot \{1 - (sleaf\,_2(t))^2 - (cleaf\,_2(t))^2\} = 0$$
(3.57)

It is necessary to satisfy one of the following conditions (i)–(iii) so that the above equation can satisfy the equation: *dx(t)/dt=0*

(i)     $sleaf\,_2(t) - cleaf\,_2(t) = 0$     (3.58)
(ii)    $sleaf\,_2(t) + cleaf\,_2(t) = 0$     (3.59)
(iii)   $1 - (sleaf\,_2(t))^2 - (cleaf\,_2(t))^2 = 0$     (3.60)

We consider the case where condition (i) is satisfied in domain (1): $2m\pi_2 \leq t < \frac{\pi_2}{2} + 2m\pi_2$.

By using Eq. (3.57) and Eq. (66) in Ref. [13], the following equation is obtained.

$$(sleaf_2(t))^2 = \sqrt{2} - 1 \quad or \quad sleaf_2(t) = \sqrt{\sqrt{2} - 1} \tag{3.61}$$

The variable $t$ is obtained as follows:

$$t = \int_0^{\sqrt{\sqrt{2}-1}} \frac{1}{\sqrt{1-u^4}} du = \frac{\pi_2}{4} (= 0.655514\cdots) \tag{3.62}$$

The minimum value of Eq. (3.52) is obtained as follows:

$$x\left(2m\pi_2 + \frac{\pi_2}{4}\right) = \sqrt{1+\sqrt{2}-1} + \sqrt{1+\sqrt{2}-1} = 2^{\frac{5}{4}} (= 2.378\cdots) \tag{3.63}$$

Next, we consider the case where condition (ii) is satisfied in domain (1): $2m\pi_2 \leq t < \frac{\pi_2}{2} + 2m\pi_2$. As shown in Fig. 19, the curve of the function: $sleaf_2(t)$ does not intersect the curve of the function: $-cleaf_2(t)$. Therefore, the variable $t$ is not satisfied with condition (ii).

Next, we consider the case where condition (iii) is satisfied in domain (1): $2m\pi_2 \leq t < \frac{\pi_2}{2} + 2m\pi_2$. The following equation is obtained by squaring Eq. (3.52).

$$\{x(t)\}^2 = 1 + (sleaf_2(t))^2 + 1 + (cleaf_2(t))^2 + 2\sqrt{1+(sleaf_2(t))^2}\sqrt{1+(cleaf_2(t))^2} \tag{3.64}$$

By using Eq. (3.59) and Eq. (66) in Ref. [13], the following equation is obtained.

$$\{x(t)\}^2 = 3 + 2\sqrt{2} \tag{3.65}$$

As shown in Fig. 18, the type(V) exact solution $x(t)$ becomes positive with respect to the arbitrary variable $t$. The following equation is then obtained:

$$x(2m\pi_2) = \sqrt{3+2\sqrt{2}} = 1 + \sqrt{2} = 2.414214 \tag{3.66}$$

Next, we consider the case where condition (i) is satisfied in domain (2): $\frac{\pi_2}{2} + 2m\pi_2 \leq t < \pi_2 + 2m\pi_2$. As shown in Fig. 19, the curve of the function: *sleaf₂(t)* does not intersect the curve of the function: -c*leaf₂(t)*. Therefore, the variable *t* satisfied with condition (i) cannot be obtained.

Next, we consider the case where condition (ii) is satisfied in domain (2): $\frac{\pi_2}{2} + 2m\pi_2 \leq t < \pi_2 + 2m\pi_2$. The minimum value of Eq. (3.52) is obtained as follows:

$$x\left(2m\pi_2 + \frac{3}{4}\pi_2\right) = \sqrt{1+\sqrt{2}-1} + \sqrt{1+\sqrt{2}-1} = 2^{\frac{5}{4}} (= 2.378\cdots) \quad (3.67)$$

The maximum value of Eq. (3.52) is obtained as follows:

$$x\left(2m\pi_2 + \frac{\pi_2}{2}\right) = 1 + \sqrt{2} = 2.414214 \quad (3.68)$$

Next, we consider the case where condition (i) is satisfied in domain (3): $\pi_2 + 2m\pi_2 \leq t < \frac{3\pi_2}{2} + 2m\pi_2$. The minimum value of Eq. (3.52) is obtained as follows:

$$x\left(2m\pi_2 + \frac{5}{4}\pi_2\right) = \sqrt{1+\sqrt{2}-1} + \sqrt{1+\sqrt{2}-1} = 2^{\frac{5}{4}} \quad (3.69)$$

Next, we consider the case where condition (ii) is satisfies in domain (3): $\pi_2 + 2m\pi_2 \leq t < \frac{3\pi_2}{2} + 2m\pi_2$. As shown in Fig. 18, the curve of the function: *sleaf₂(t)* does not intersect the curve of the function: -c*leaf₂(t)*. Therefore, the variable *t* satisfied with condition (i) cannot be obtained.

Next, we consider the case where condition (iii) is satisfied in domain (3): $\pi_2 + 2m\pi_2 \leq t < \frac{3\pi_2}{2} + 2m\pi_2$. The maximum value of Eq. (3.52) is obtained as follows:

$$x(2m\pi_2 + \pi_2) = 1 + \sqrt{2} \quad (3.70)$$

Next, we consider the case where condition (ii) is satisfied in domain (4): $\frac{3\pi_2}{2}+2m\pi_2 \leq t < 2\pi_2+2m\pi_2$. As shown in Fig. 19, the curve of the function: $sleaf_2(t)$ does not intersect the curve of the function: $cleaf_2(t)$. Therefore, the variable $t$ satisfied with condition (i) cannot be obtained.

Next, we consider the case where condition (ii) is satisfied in domain (4): $\frac{3\pi_2}{2}+2m\pi_2 \leq t < 2\pi_2+2m\pi_2$. The minimum value of Eq. (3.52) is obtained as follows:

$$x\left(2m\pi_2+\frac{7}{4}\pi_2\right)=\sqrt{1+\sqrt{2}-1}+\sqrt{1+\sqrt{2}-1}=2^{\frac{5}{4}} \quad (3.71)$$

Next, we consider the case where condition (iii) is satisfied in domain (4): $\frac{3\pi_2}{2}+2m\pi_2 \leq t < 2\pi_2+2m\pi_2$. The maximum value of Eq. (3.52) is obtained as follows:

$$x\left(2m\pi_2+\frac{3\pi_2}{2}\right)=1+\sqrt{2} \quad (3.72)$$

As stated above, the range of the variable $x(t)$ is as follows:

$$2^{\frac{5}{4}} \leq x(t) \leq 1+\sqrt{2} \quad (3.73)$$

The center of the displacement and amplitude becomes $x(t)=\frac{1+\sqrt{2}+2^{\frac{5}{4}}}{2}$ and $x(t)=\frac{1+\sqrt{2}-2^{\frac{5}{4}}}{2}$, respectively.

Next, the variables are set to $\Phi=0$ and $\omega=1$. The amplitude $A$ is varied. Under these conditions, the wave obtained by the type (Ⅴ) exact solution is shown in Fig. 20. The range of the displacement $x(t)$ can be obtained by the following inequality:

$$2^{\frac{5}{4}}A \leq x(t) \leq (1+\sqrt{2})A \quad A>0 \quad (3.74)$$

$$(1+\sqrt{2})A \leq x(t) \leq 2^{\frac{5}{4}}A \quad A<0 \tag{3.75}$$

The center of the displacement *x(t)* is obtained as follows:

$$(Center\ of\ displacement) = \frac{2^{\frac{5}{4}}+1+\sqrt{2}}{2}A \tag{3.76}$$

The amplitude is obtained as follows:

$$(Amplitude) = (1+\sqrt{2})A - \frac{2^{\frac{5}{4}}+1+\sqrt{2}}{2}|A| = \frac{1+\sqrt{2}-2^{\frac{5}{4}}}{2}|A| \tag{3.77}$$

Next, the variables are set to *Φ=0* and *A=1*. The angular frequency *ω* is varied. The wave obtained by the type (Ⅴ) exact solution is shown in Fig. 21. The period of the waves is varied according to the absolute value *ω*. As the absolute value *ω* increases, the period decreases. In contrast, as the absolute value *ω* decreases, the period increases. In case *ω=±1*, the period becomes the constant $\pi_2/2$. In case *ω=±2*, the period becomes the constant $\pi_2/4$. In case *ω=±3*, the period becomes the constant $\pi_2/6$. By using the parameter *ω*, the period *T* is obtained as follows:

$$T = \frac{\pi_2}{2\omega} \tag{3.78}$$

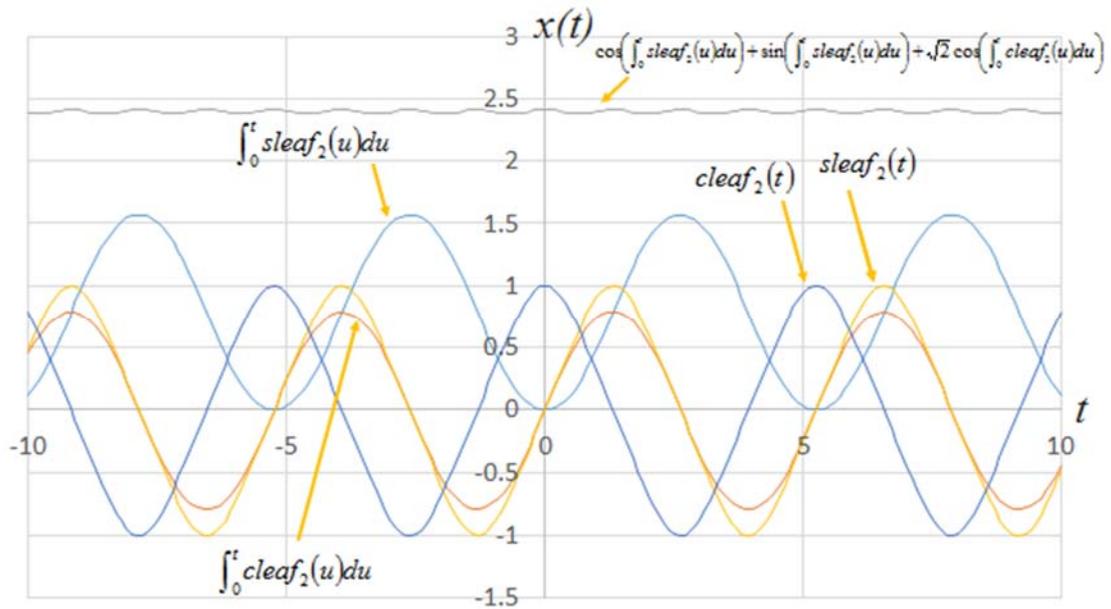

Fig. 17 Waves of the exact solution obtained by Type(V); the leaf function: $sleaf_2(t)$ and the function: $\int_0^t sleaf_2(u)du$

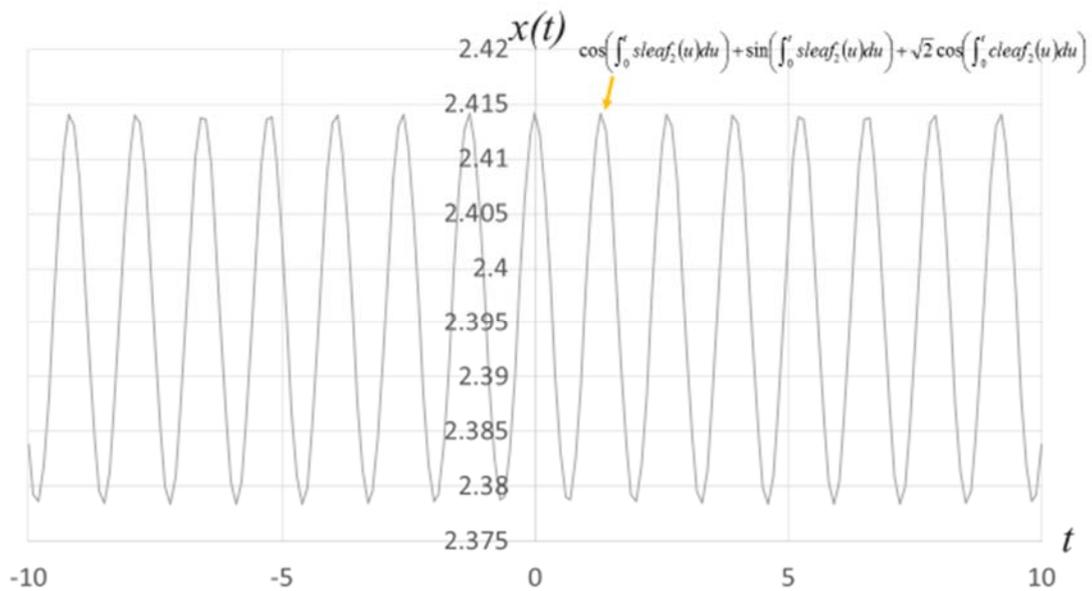

Fig. 18 Enlargement of Fig. 17

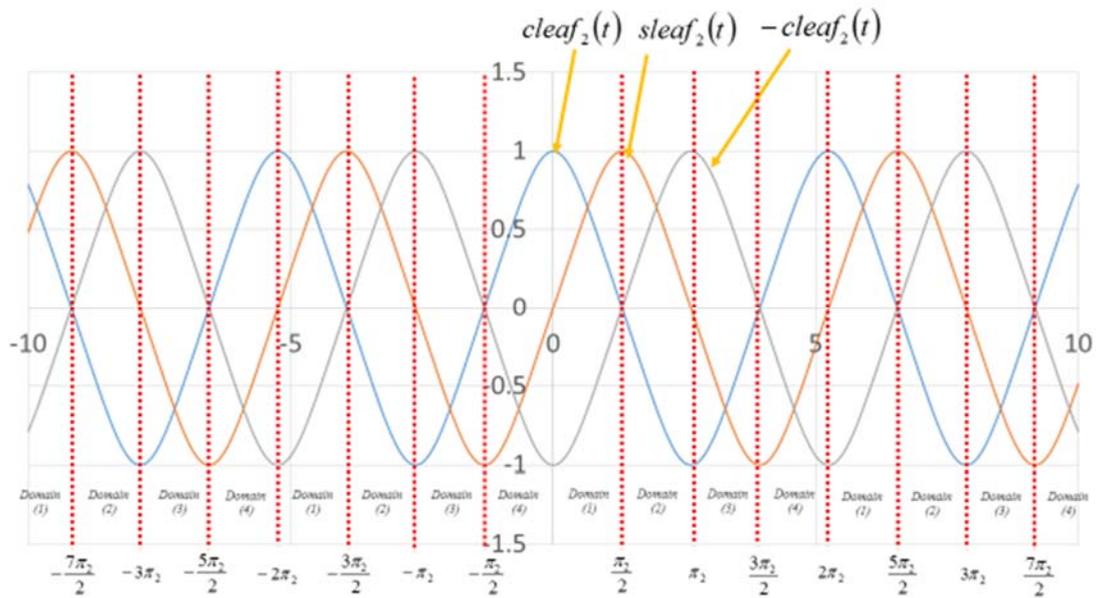

Fig.19 Waves obtained by the function: *cleaf₂(t)*, *sleaf₂(t)* and -*cleaf₂(t)*

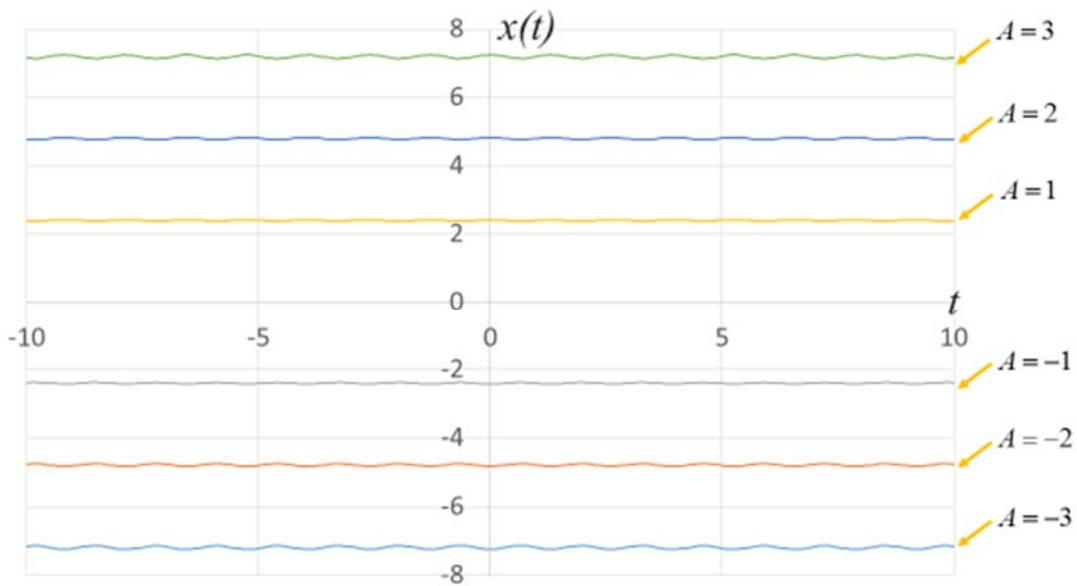

Fig. 20 Wave obtained by Type(V) exact solution with respect to the variation in the amplitude $A$

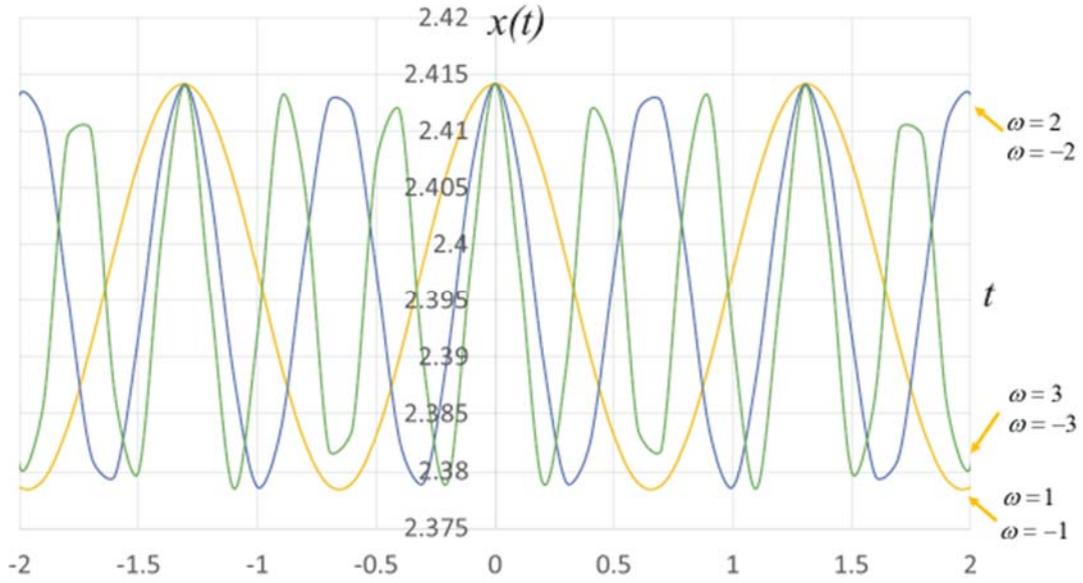

Fig. 21 Wave obtained by Type (V) exact solution with respect to the variation in the angular frequency $\omega$

### 3.6 Numerical Results of Exact Solution of Type (VI)

In the case $A=1$, $\omega=1$, $\Phi=0$ in Eq. (2.22), the wave of the exact solution of type (VI), the wave of the leaf function: $sleaf_2(t)$ and $cleaf_2(t)$, the wave of the function: $\int_0^t sleaf_2(u)du$ and $\int_0^t cleaf_2(u)du$ are shown in Fig. 22. The horizontal and vertical axes represent time and displacement $x(t)$, respectively. The numerical data of the Type (VI) solution can be obtained by using the table 3 and the table 5. To discuss the range of the type (VI) exact solution $x(t)$, the following equation is transformed by using the leaf function: $sleaf_2(t)$ and $cleaf_2(t)$:

$$x(t) = \sqrt{1+(sleaf_2(t))^2} - \sqrt{1+(cleaf_2(t))^2} \quad (3.79)$$

The above equation can be obtained by the same operation from Eq. (3.46) to Eq. (3.51). To discuss the range of variable $x(t)$ in the type (VI) exact solution, the first order differential is obtained. The sign of the first order differential with respect to the leaf function $sleaf_2(t)$ and $cleaf_2(t)$ depends on the domain $t$ of the variable $x(t)$. The first order differential of Eq. (3.79) is discussed by division into four domains (1)–(4).

Domain (1): $2m\pi_2 \leq t < \dfrac{\pi_2}{2} + 2m\pi_2$ (m:integer)

$$\dfrac{dx(t)}{dt} = sleaf_2(t)\sqrt{1-(sleaf_2(t))^2} + cleaf_2(t)\sqrt{1-(cleaf_2(t))^2} \tag{3.80}$$

Domain (2): $\dfrac{\pi_2}{2} + 2m\pi_2 \leq t < \pi_2 + 2m\pi_2$ (m:integer)

$$\dfrac{dx(t)}{dt} = -sleaf_2(t)\sqrt{1-(sleaf_2(t))^2} + cleaf_2(t)\sqrt{1-(cleaf_2(t))^2} \tag{3.81}$$

Domain (3): $\pi_2 + 2m\pi_2 \leq t < \dfrac{3\pi_2}{2} + 2m\pi_2$ (m:integer)

$$\dfrac{dx(t)}{dt} = -sleaf_2(t)\sqrt{1-(sleaf_2(t))^2} - cleaf_2(t)\sqrt{1-(cleaf_2(t))^2} \tag{3.82}$$

Domain (4): $\dfrac{3\pi_2}{2} + 2m\pi_2 \leq t < 2\pi_2 + 2m\pi_2$ (m:integer)

$$\dfrac{dx(t)}{dt} = sleaf_2(t)\sqrt{1-(sleaf_2(t))^2} - cleaf_2(t)\sqrt{1-(cleaf_2(t))^2} \tag{3.83}$$

The extreme value of the type(VI) solution is obtained by the equation: $dx(t)/dt=0$ with respect to Eqs. (3.80)–(3.83). The same equation of Eq. (3.57) is derived from the equation: $dx(t)/dt=0$ with respect to Eqs. (3.80)–(3.83). Therefore, in order to satisfy the equation: $dx(t)/dt=0$, it is necessary to satisfy any of conditions (i)–(iii), as shown in Eqs. (3.58)–(3.60).

Next, we consider the case where condition (i) is satisfied in domain (1): $2m\pi_2 \leq t < \dfrac{\pi_2}{2} + 2m\pi_2$. As shown in Fig. 19, the curve of the function: $sleaf_2(t)$ does not intersect the curve of the function: $cleaf_2(t)$. Therefore, the variable $t$ satisfied with condition (i) cannot be obtained. The extremal value of Eq. (3.80) cannot be obtained under condition (i).

Next, we consider the case where condition (ii) is satisfied in domain (1): $2m\pi_2 \leq t < \dfrac{\pi_2}{2} + 2m\pi_2$. As shown in Fig. 19, the curve of the function: $sleaf_2(t)$ does not intersect the curve of the function: $-cleaf_2(t)$. Therefore, the variable $t$ satisfied with condition (i) cannot be obtained. The extremal value of Eq. (3.80) cannot be obtained under condition (ii).

Next, we consider the case where condition (iii) satisfies in domain (1): $2m\pi_2 \le t < \frac{\pi_2}{2} + 2m\pi_2$. The following equation is obtained by squaring both sides of Eq. (3.79).

$$\{x(t)\}^2 = 1 + (sleaf_2(t))^2 + 1 + (cleaf_2(t))^2 - 2\sqrt{1 + (sleaf_2(t))^2}\sqrt{1 + (cleaf_2(t))^2} \tag{3.84}$$

Using Eq. (3.60) and Eq. (66) in Ref. [13], the following equation is obtained.

$$\{x(t)\}^2 = 3 - 2\sqrt{2} \tag{3.85}$$

The following equation is obtained by the above equation.

$$x(2m\pi_2) = -\sqrt{3 - 2\sqrt{2}} = -(\sqrt{2} - 1) = -0.414214\cdots \quad \text{(m: integer)} \tag{3.86}$$

Next, we consider the case where condition (i) is satisfied in domain (2): $\frac{\pi_2}{2} + 2m\pi_2 \le t < \pi_2 + 2m\pi_2$. As shown in Fig. 19, the curve of the function: $sleaf_2(t)$ does not intersect the curve of the function: $cleaf_2(t)$. Therefore, the variable $t$ satisfied with condition (ii) cannot be obtained. The extremal value of Eq. (3.81) cannot be obtained under condition (i).

Next, we consider the case where condition (ii) is satisfied in domain (2): $\frac{\pi_2}{2} + 2m\pi_2 \le t < \pi_2 + 2m\pi_2$. As shown in Fig. 19, the intersection points between the function: $sleaf_2(t)$ and the function: $-cleaf_2(t)$ are shown in Eq. (3.61). The intersection point described in Eq. (3.61) is not an extremal value, so Eq. (3.81) does not become $0$.

Next, we consider the case where condition (iii) is satisfied in domain (2): $\frac{\pi_2}{2} + 2m\pi_2 \le t < \pi_2 + 2m\pi_2$. The following equation is obtained by squaring both sides of Eq. (3.79).

$$x\left(2m\pi_2 + \frac{\pi_2}{2}\right) = \sqrt{3 - 2\sqrt{2}} = \sqrt{2} - 1 = 0.414214\cdots \quad \text{(m: integer)} \tag{3.87}$$

Next, we consider the case where condition (i) is satisfied in domain (3):

$\pi_2 + 2m\pi_2 \leq t < \frac{3\pi_2}{2} + 2m\pi_2$. As shown in Fig. 19, the intersection points between the function: *sleaf₂(t)* and the function: *-cleaf₂(t)* are given as follows:

$$sleaf_2\left(2m\pi_2 + \frac{3\pi_2}{4}\right) = cleaf_2\left(2m\pi_2 + \frac{3\pi_2}{4}\right) = -\sqrt{\sqrt{2}-1} \quad \text{(m: integer)} \tag{3.88}$$

The above equation is not satisfied with *dx/dt=0* with respect to Eq. (3.82). The extremal value of Eq. (3.82) cannot be obtained under condition (i).

Next, we consider the case where condition (ii) is satisfied in domain (3): $\pi_2 + 2m\pi_2 \leq t < \frac{3\pi_2}{2} + 2m\pi_2$. As shown in Fig. 19, the curve of the function: *sleaf₂(t)* does not intersect the curve of the function: *-cleaf₂(t)*. Therefore, the variable *t* satisfied with condition (ii) cannot be obtained.

Next, we consider the case where condition (iii) is satisfied in domain (3): $\pi_2 + 2m\pi_2 \leq t < \frac{3\pi_2}{2} + 2m\pi_2$. The variable *x(t)* is obtained by squaring both sides of Eq. (3.79).

$$x(2m\pi_2 + \pi_2) = -\left(\sqrt{2}-1\right) \quad \text{(m: integer)} \tag{3.89}$$

Next, we consider the case where condition (i) is satisfied in domain (4): $\frac{3\pi_2}{2} + 2m\pi_2 \leq t < 2\pi_2 + 2m\pi_2$. As shown in Fig. 19, the curve of the function: *sleaf₂(t)* does not intersect the curve of the function: *cleaf₂(t)*. Therefore, the variable *t* satisfied with condition (i) cannot be obtained. The extremal value of Eq. (3.83) cannot be obtained under condition (i).

Next, we consider the case where condition (ii) is satisfied in domain (4): $\frac{3\pi_2}{2} + 2m\pi_2 \leq t < 2\pi_2 + 2m\pi_2$. As shown in Fig. 19, the intersection points between the function: *sleaf₂(t)* and the function: *- cleaf₂(t)* are given as follows:

$$sleaf_2\left(2m\pi_2 + \frac{5\pi_2}{4}\right) = cleaf_2\left(2m\pi_2 + \frac{5\pi_2}{4}\right) = -\sqrt{\sqrt{2}-1} \quad \text{(m: integer)} \tag{3.90}$$

The intersection point described in the above equation is not an extremal value, so Eq. (3.83) does not become *0*. Next, we consider the case where condition (iii) is satisfied in domain (4): $\frac{3\pi_2}{2} + 2m\pi_2 \leq t < 2\pi_2 + 2m\pi_2$. The variable $x(t)$ is obtained by squaring both sides of Eq. (3.79).

$$x\left(2m\pi_2 + \frac{3\pi_2}{2}\right) = \sqrt{2} - 1 \quad (m: integer) \tag{3.91}$$

As stated above, the range of the variable $x(t)$ is as follows:

$$-\left(\sqrt{2} - 1\right) \leq x(t) \leq \sqrt{2} - 1 \tag{3.92}$$

The amplitude becomes $\sqrt{2} - 1$. The variables are set to $\Phi=0$ and $\omega=1$. The amplitude A is varied. Under these conditions, the wave obtained by the type (VI) exact solution is shown in Figs. 23 and 24. The range of the displacement $x(t)$ can be obtained by the following inequality:

$$-\left(\sqrt{2} - 1\right)A \leq x(t) \leq \left(\sqrt{2} - 1\right)A \tag{3.93}$$

The center of the displacement $x(t)$ and amplitude become $0.0$ and $\left(\sqrt{2} - 1\right)A$, respectively.
Next, the variables are set to $\Phi=0$ and $A=1$. The angular frequency $\omega$ is varied. The wave obtained by the type (VI) exact solution is shown in Fig. 25. The period of the waves is varied according to the absolute value $\omega$. As the absolute value $\omega$ increases, the period decreases. In contrast, as the absolute value $\omega$ decreases, the period increases. In case $\omega=\pm 1$, the period becomes the constant $\pi_2$. In case $\omega=\pm 2$, the period becomes the constant $\pi_2/2$. In case $\omega=\pm 3$, the period becomes the constant $\pi_2/3$. By using the parameter $\omega$, the period $T$ is obtained as follows:

$$T = \frac{\pi_2}{|\omega|} \tag{3.94}$$

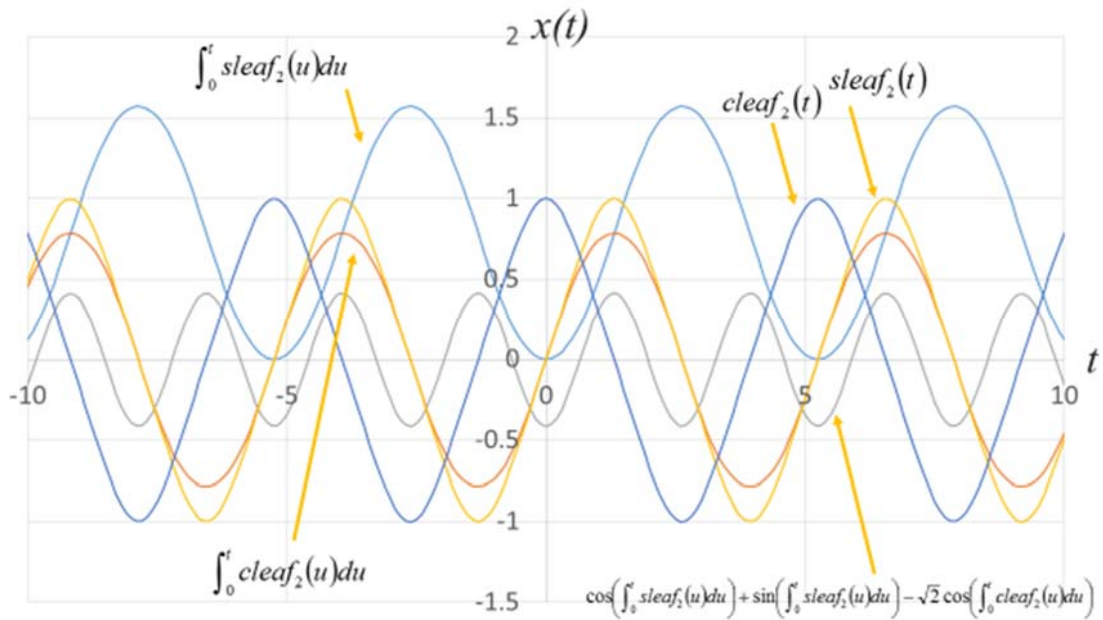

Fig. 22 Waves of the exact solution obtained by Type(VI); the leaf function: $sleaf_2(t)$ and the function: $\int_0^t sleaf_2(u)du$

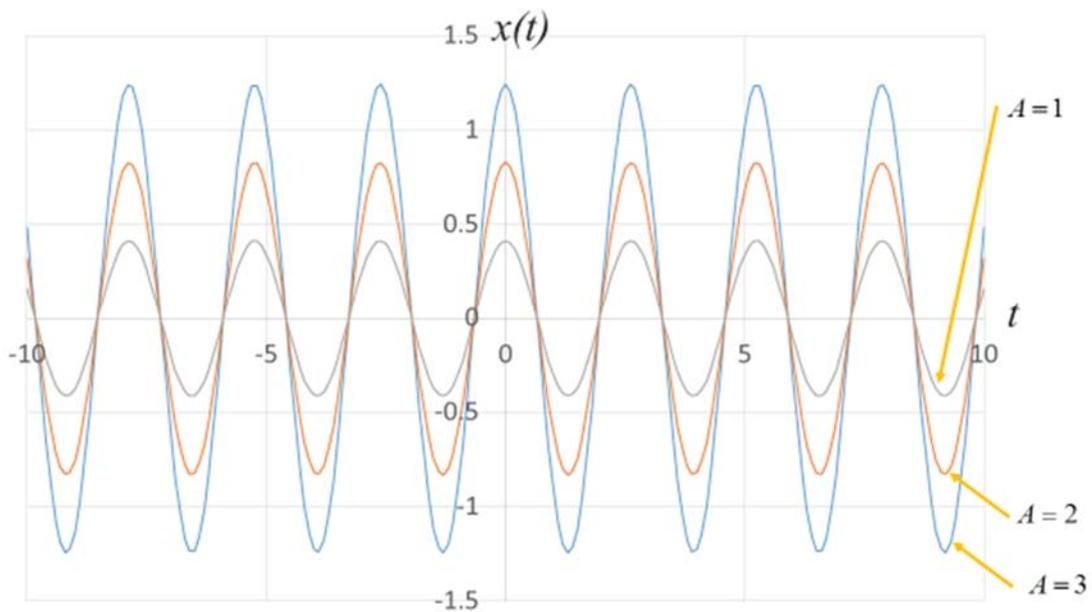

Fig. 23 Wave obtained by Type(VI) exact solution with respect to the variation in the amplitude $A$ ($A=1,2,3$)

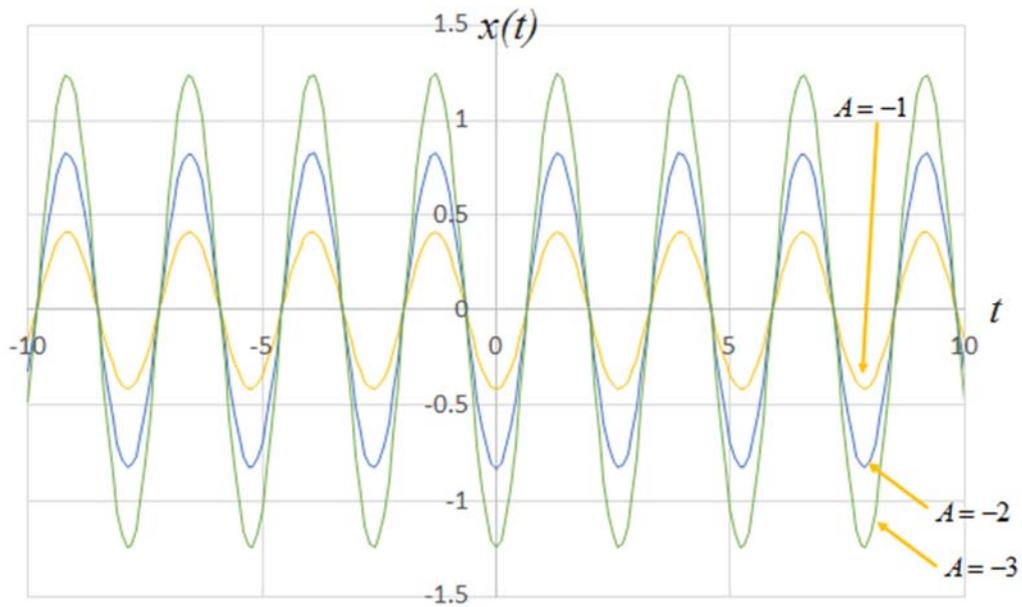

Fig. 24 Wave obtained by Type(VI) exact solution with respect to the variation in the amplitude $A$ ($A$=-1,-2,-3)

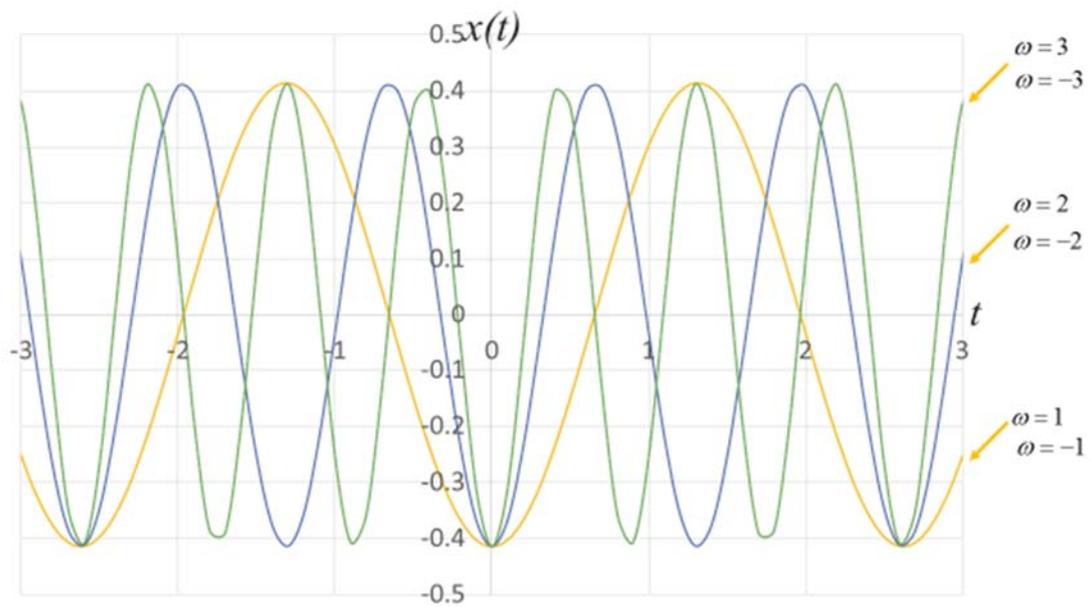

Fig. 25 Wave obtained by Type (VI) exact solution with respect to the variation in the angular frequency $\omega$

## 3.7 Numerical Results of Exact Solution of Type( Ⅶ)

In the case $A=1$, $\omega=1$, $\Phi=0$ in Eq. (2.26), the wave of the exact solution of type (Ⅶ)：$sleaf_2(t) \cdot cleaf_2(t)$, the wave of the leaf function: $sleaf_2(t)$ and $cleaf_2(t)$ are shown in Fig. 26. The horizontal and vertical axes represent time and displacement $x(t)$, respectively. The numerical data of the Type (Ⅶ) solution can be obtained by using the table 1. Next, the variables are set to $\Phi=0$ and $\omega=1$. The amplitude $A$ is varied. Under these conditions, the wave obtained by the type (Ⅶ) exact solution is shown in Figs. 27-28. In order to obtain the extreme values, the first order differential is derived as follows:

$$\frac{dx(t)}{dt} = A \cdot cleaf_2(t) \cdot \sqrt{1-(sleaf_2(t))^4} - A \cdot sleaf_2(t)\sqrt{1-(cleaf_2(t))^4} = 0 \tag{3.95}$$

The following equation is derived as follows:

$$cleaf_2(t) = \pm sleaf_2(t) \tag{3.96}$$

By using the above equation and Eq. (66) in Ref. [2], the function $sleaf_2(t)$ and $cleaf_2(t)$ is obtained as follows:

$$cleaf_2(t) = sleaf_2(t) = \pm\sqrt{-1+\sqrt{2}} \tag{3.97}$$

From the above, the minimum of the variable $x(t)$ is obtained as follows:

$$x\left(\frac{\pi_2}{4}(4m+1)\right) = A\left(\sqrt{-1+\sqrt{2}}\right)^2 = A \cdot \left(-1+\sqrt{2}\right) \quad (m=0, \pm 1, \pm 2 \cdots) \tag{3.98}$$

The maximum of the variable $x(t)$ is obtained as follows:

$$x\left(\frac{\pi_2}{4}(4m+3)\right) = -A\left(\sqrt{-1+\sqrt{2}}\right)^2 = -A \cdot \left(-1+\sqrt{2}\right) \quad (m=0, \pm 1, \pm 2 \cdots) \tag{3.99}$$

The range of displacement of x(t) can be obtained by following inequlity:

$$-\left(-1+\sqrt{2}\right)A \leq x(t) \leq \left(-1+\sqrt{2}\right)A \tag{3.100}$$

The angular frequency $\omega$ is varied. The wave obtained by the type (VII) exact solution is shown in Figs. 29 and 30. As shown in Figs. 29 and 30, the period $T$ is obtained from the above results.

$$T = \frac{\pi_2}{|\omega|} \tag{3.101}$$

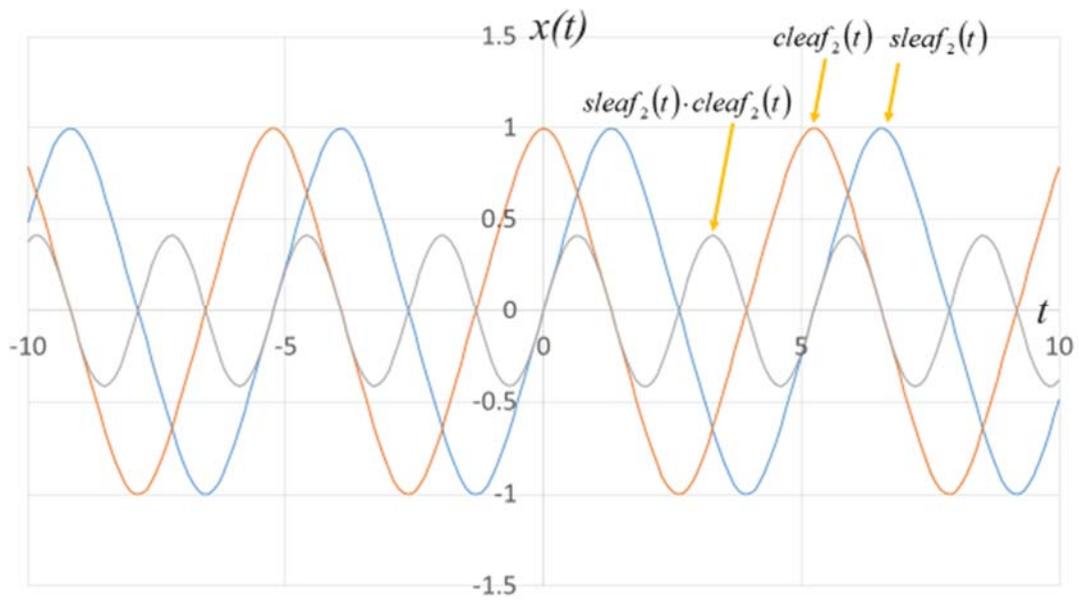

Fig. 26 Waves of the exact solution obtained by Type(VII); the leaf function: $sleaf_2(t)$ and $cleaf_2(t)$

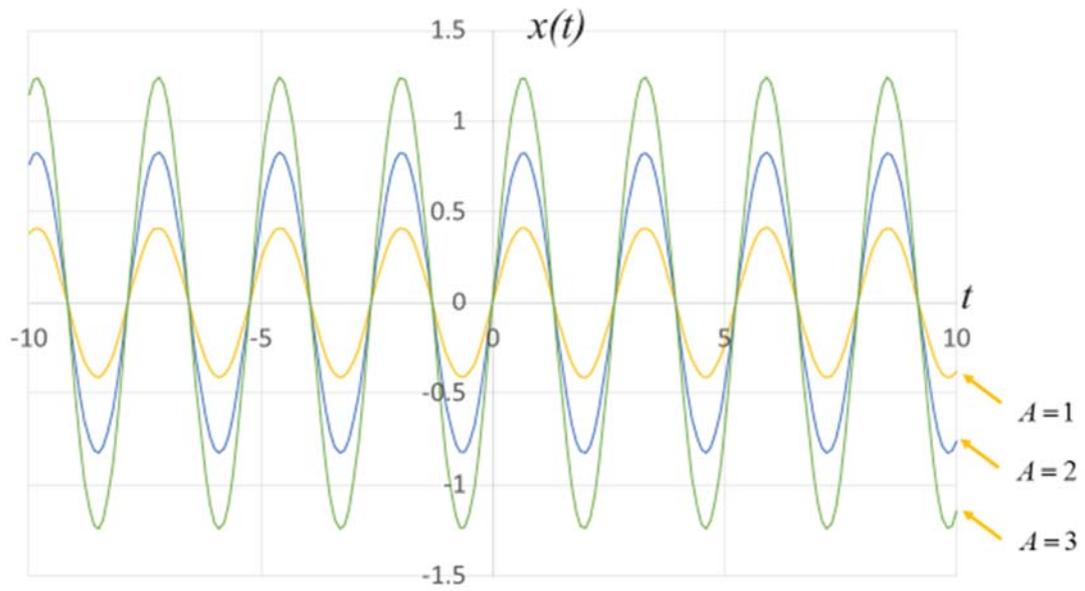

Fig. 27 Wave obtained by Type(VII) exact solution with respect to the variation in the amplitude $A(A=1,2,3)$

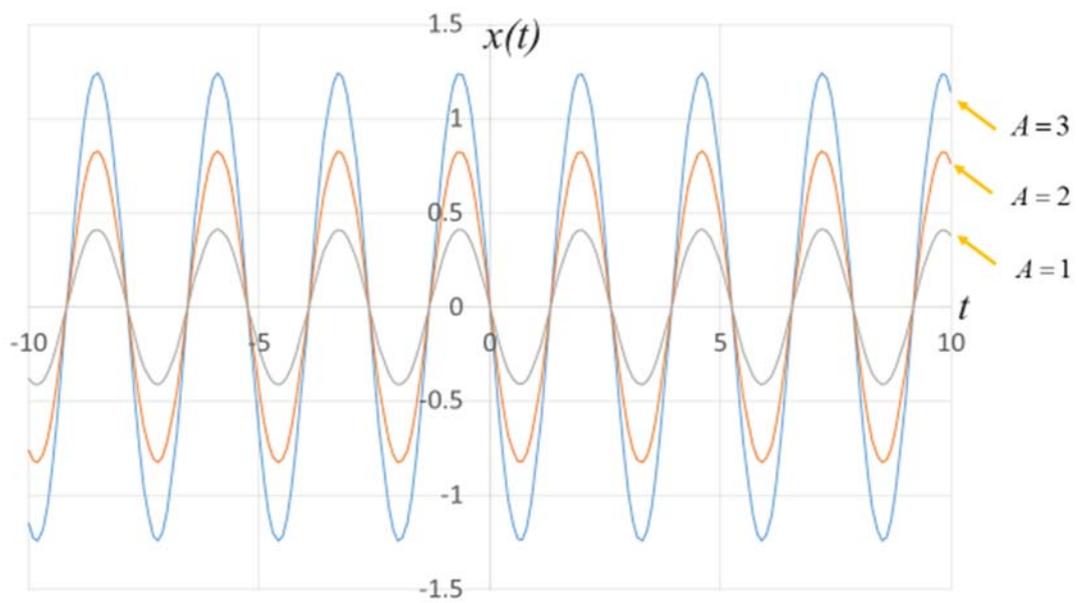

Fig. 28 Wave obtained by Type(VII) exact solution with respect to the variation in the amplitude $A(A=-1, -2, -3)$

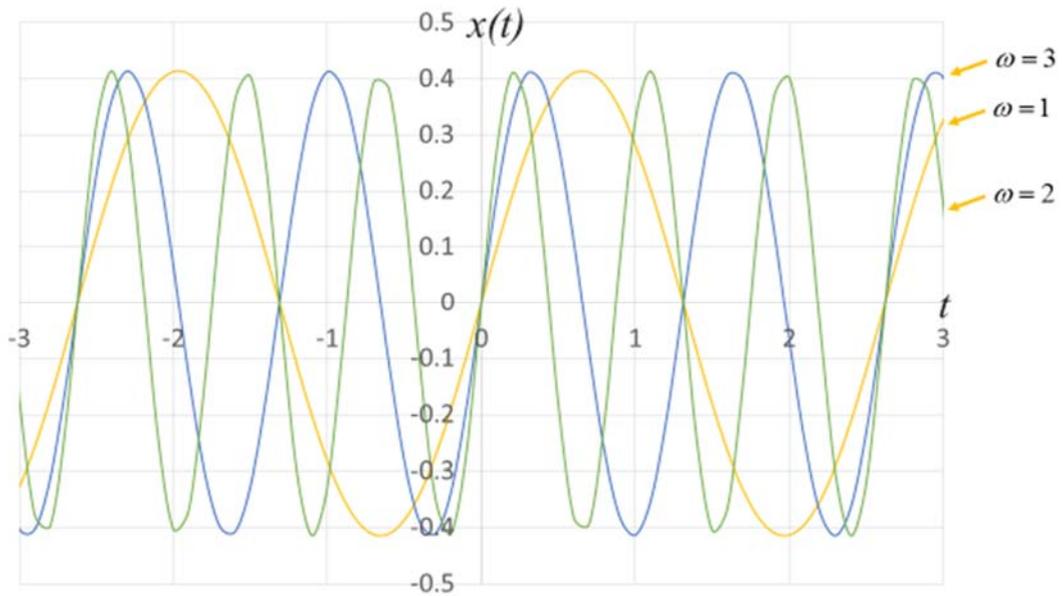

Fig. 29 Wave obtained by Type (VII) exact solution with respect to the variation in the angular frequency $\omega$ ($\omega=1, 2, 3$)

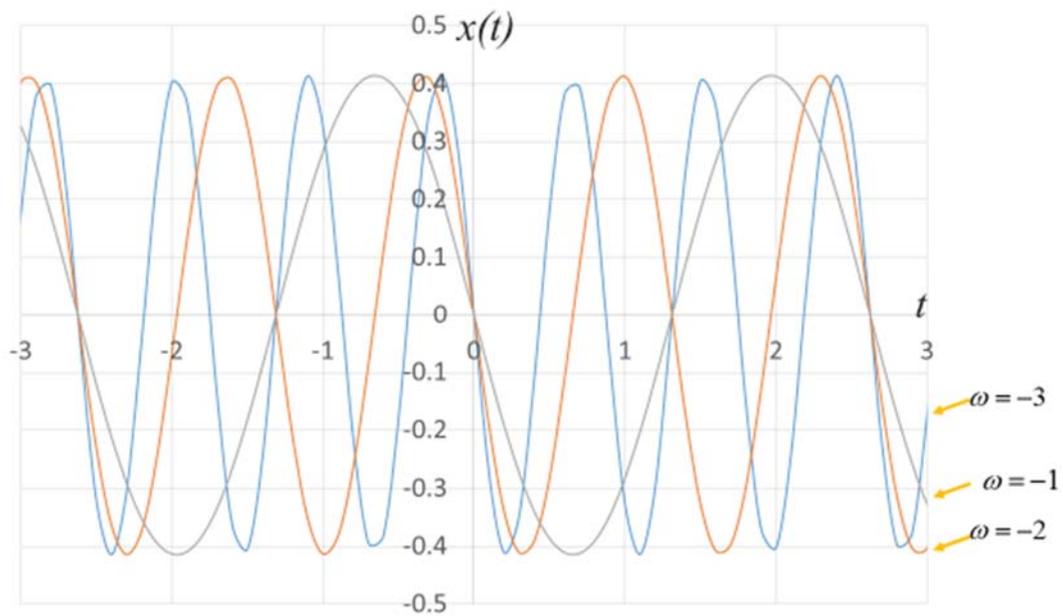

Fig. 30 Wave obtained by Type (VII) exact solution with respect to the variation in the angular frequency $\omega$ ($\omega=-1, -2, -3$)

*4. Conclusion*

  Leaf functions are functions in which the second order differential of the leaf function is equal to *-n* times the function raised to the *(2n-1)* power of the leaf function. By using leaf functions, the exact solution of the cubic duffing equation can be derived under certain conditions. With the exact solution, the wave obtained by the exact solution is visualized in a graph. The conclusions are summarized as follows:

  ・Through leaf functions, seven types of exact solutions can be derived from the cubic duffing equation.
  ・The seven types of exact solutions have two parameters, namely the angular frequency $\omega$ and the amplitude *A*. The parameters *A* and $\omega$ indicate the characteristics of the wave. The coefficients of the terms *x* and $x^3$ in the cubic duffing equation can be described by both the wave amplitude *A* and wave frequency parameters $\omega$ in the leaf functions. The amplitude *A* of the wave becomes constant, even if these coefficients vary according to the variation of the parameter $\omega$. In contrast, the amplitude $\omega$ of the wave becomes constant, even if these coefficients vary according to the variation in the parameter *A*. Since parameters *A* and ω do not affect the characteristics of the wave, they are independent variables in the ordinary differential equation.
  ・As the variable representing the amplitude increases, the height of the wave also increases. In contrast, as the variable representing the amplitude decreases, the height of the wave also decreases. As the variable $\omega$ representing the period increases, the period of the waves decreases. In contrast, as the variable $\omega$ representing the period decreases, the period of the wave increases. The waveform obtained by the nonlinear spring can be controlled by adjusting these variables. Several new waveforms satisfying the cubic duffing equation can be constructed by combining both trigonometric functions and leaf functions.

*References*

[1] K. I. Brennan et al., The Duffing Equation: Nonlinear Oscillators and their Behaviour, *Wiley*, 2011.
[2] C. Livija, Ninety years of Duffing's equation, *Theoretical and Applied Mechanics*, **40**(2013), 49-63.
[3] A. Beléndez et al., "Analytical Approximate Solutions for the Cubic-Quintic Duffing Oscillator in Terms of Elementary Functions," Journal of Applied Mathematics, Article

ID 286290, (2012).

[4] A. Beléndez et al., An accurate closed-form approximate solution for the quintic Duffing oscillator equation, Mathematical and Computer Modelling, **52**(2010), 637-641

[5] U. Starossek, Exact analytical solutions for forced undamped Duffing oscillator, International Journal of Non-Linear Mechanics, **85**(2016), 197-206.

[6] H. Alvaro et al., Exact solution to duffing equation and the pendulum equation, Applied Mathematical Sciences, **8**(2014), 8781-8789

[7] S. Salas et al., Exact solutions to cubic duffing equation for a nonlinear electrical circuit, Visión Electrónica: algo más que un estado sólido, **8**(2014), p. 46-53.

[8] L. Cveticanin, The approximate solving methods for the cubic duffing equation based on the Jacobi elliptic functions, International Journal of Nonlinear Sciences and Numerical Simulation, **10**(2009), 1491-1516.

[9] V. Marinca et al., "Explicit and exact solutions to cubic duffing and double-well duffing equations", Mathematical and Computer Modelling, **53**(2011), 604-609.

[10] B. T. Nohara et al., Solution of the duffing equation with a higher order nonlinear term, Theoretical and Applied Mechanics Japan, **59**(2011), 133-141.

[11] S. B. Yamgoué et al., Analytical solutions of undamped and autonomous cubic-quintic duffing equation, American Journal of Physics and Applications. **3**(2015), pp. 159-165

[12] K. Shinohara, Special function: Leaf function $r=sleaf_n(l)$ (First report), Bulletin of Daido University. (2015), 51;23–38.
URL: http://id.nii.ac.jp/1277/00000049/

[13] K. Shinohara, Special function: Leaf function $r=cleaf_n(l)$ (Second report), Bulletin of Daido University, (2015), 51;39–68.
URL: http://id.nii.ac.jp/1277/00000050/

*Appendix I*

The type(I) exact solution is satisfied with the cubic duffing equation in Eq. (2.1). Eq. (2.2) is given as follows:

$$x(t) = A\cos\left(\int_0^{\omega t+\phi} cleaf_2(u)du\right)$$

( I .1)

With the above equation, the first order differential is obtained as follows:

$$\frac{dx(t)}{dt} = -A\sin\left(\int_0^{\omega t+\phi} cleaf_2(u)du\right)\cdot \omega \cdot cleaf_2(\omega t+\phi) \tag{I.2}$$

With the above equation, the second order differential is obtained as follows:

$$\frac{d^2x(t)}{dt^2} = -A\cos\left(\int_0^{\omega t+\phi} cleaf_2(u)du\right)\cdot \omega^2 \cdot (cleaf_2(\omega t+\phi))^2$$
$$- A\sin\left(\int_0^{\omega t+\phi} cleaf_2(u)du\right)\cdot \omega^2 \cdot \left(-\sqrt{1-(cleaf_2(\omega t+\phi))^4}\right) \tag{I.3}$$

By using Eq. (63) in Ref. [13], the following equation is obtained:

$$(cleaf_2(\omega t+\phi))^2 = \cos\left(2\int_0^{\omega t+\phi} cleaf_2(u)du\right) \tag{I.4}$$

The above equation is transformed as follows:

$$\sqrt{1-(cleaf_2(\omega t+\phi))^4} = \sqrt{1-\left(\cos\left(2\int_0^{\omega t+\phi} cleaf_2(u)du\right)\right)^2} = \sin\left(2\int_0^{\omega t+\phi} cleaf_2(u)du\right) \tag{I.5}$$

The following equation is obtained from the above equation:
Substituting the above two equations into Eq. (I.3), the following equation is obtained:

$$\frac{d^2x(t)}{dt^2} = -A\omega^2 \cos\left(\int_0^{\omega t+\phi} cleaf_2(u)du\right)\cdot \cos\left(2\int_0^{\omega t+\phi} cleaf_2(u)du\right)$$
$$+ A\omega^2 \sin\left(\int_0^{\omega t+\phi} cleaf_2(u)du\right)\cdot \sin\left(2\int_0^{\omega t+\phi} cleaf_2(u)du\right)$$
$$= -A\omega^2 \cos\left(\int_0^{\omega t+\phi} cleaf_2(u)du + 2\int_0^{\omega t+\phi} cleaf_2(u)du\right)$$
$$= -A\omega^2 \cos\left(3\int_0^{\omega t+\phi} cleaf_2(u)du\right)$$
$$= 3\omega^2 A\cos\left(\int_0^{\omega t+\phi} cleaf_2(u)du\right) - 4\left(\frac{\omega}{A}\right)^2\left\{A\cos\left(\int_0^{\omega t+\phi} cleaf_2(u)du\right)\right\}^3 \tag{I.6}$$

Substituting Eq. (2.2) into the above equation, the following equation is obtained:

$$\frac{d^2 x(t)}{dt^2} = 3\omega^2 x(t) - 4\left(\frac{\omega}{A}\right)^2 \{x(t)\}^3 \tag{I.7}$$

Eq. (2.3) is obtained from the above equation.

$$\frac{d^2 x(t)}{dt^2} - 3\omega^2 x(t) + 4\left(\frac{\omega}{A}\right)^2 x(t)^3 = 0 \tag{I.8}$$

The coefficients $\alpha$ and $\beta$ in Eq. (2.1) are expressed as follows:

$$\alpha = -3\omega^2 \tag{I.9}$$

$$\beta = 4\left(\frac{\omega}{A}\right)^2 \tag{I.10}$$

### Appendix II

The type(II) exact solution is satisfied with the cubic duffing equation in Eq. (2.1). Eq. (2.6) is given as follows:

$$x(t) = A \sin\left(\int_0^{\omega t + \phi} cleaf_2(u) du\right) \tag{II.1}$$

With the above equation, the first order differential is obtained as follows:

$$\frac{dx(t)}{dt} = A \cos\left(\int_0^{\omega t + \phi} cleaf_2(u) du\right) \cdot \omega \cdot cleaf_2(\omega t + \phi) \tag{II.2}$$

With the above equation, the second order differential is obtained as follows:

$$\frac{d^2x(t)}{dt^2} = -A\sin\left(\int_0^{\omega t+\phi} cleaf_2(u)du\right)\cdot\omega^2\cdot(cleaf_2(\omega t+\phi))^2$$
$$+ A\cos\left(\int_0^{\omega t+\phi} cleaf_2(u)du\right)\cdot\omega^2\cdot\left(-\sqrt{1-(cleaf_2(\omega t+\phi))^4}\right)$$
(Ⅱ.3)

Substituting Eq. (Ⅰ.4) and (Ⅰ.5) into the above equations, the following equation is obtained:

$$\frac{d^2x(t)}{dt^2} = -A\omega^2\sin\left(\int_0^{\omega t+\phi} cleaf_2(u)du\right)\cdot\cos\left(2\int_0^{\omega t+\phi} cleaf_2(u)du\right)$$
$$- A\omega^2\cos\left(\int_0^{\omega t+\phi} cleaf_2(u)du\right)\cdot\sin\left(2\int_0^{\omega t+\phi} cleaf_2(u)du\right)$$
$$= -A\omega^2\sin\left(\int_0^{\omega t+\phi} cleaf_2(u)du + 2\int_0^{\omega t+\phi} cleaf_2(u)du\right)$$
$$= -A\omega^2\sin\left(3\int_0^{\omega t+\phi} cleaf_2(u)du\right)$$
$$= -3\omega^2 A\sin\left(\int_0^{\omega t+\phi} cleaf_2(u)du\right) + 4\left(\frac{\omega}{A}\right)^2\left\{A\sin\left(\int_0^{\omega t+\phi} cleaf_2(u)du\right)\right\}^3$$
(Ⅱ.4)

Substituting Eq. (Ⅱ.1) into the above equations, the following equation is obtained:

$$\frac{d^2x(t)}{dt^2} = -3\omega^2 x(t) + 4\left(\frac{\omega}{A}\right)^2\{x(t)\}^3$$
(Ⅱ.5)

Eq. (2.7) is obtained from the above equation.

$$\frac{d^2x(t)}{dt^2} + 3\omega^2 x(t) - 4\left(\frac{\omega}{A}\right)^2 x(t)^3 = 0$$
(Ⅱ.6)

The coefficients $\alpha$ and $\beta$ in Eq. (2.1) are expressed as follows:

$$\alpha = 3\omega^2$$
(Ⅱ.7)

$$\beta = -4\left(\frac{\omega}{A}\right)^2$$
(Ⅱ.8)

*Appendix III*

The type(III) exact solution is satisfied with the cubic duffing equation in Eq. (2.1). Eq. (2.10) is given as follows:

$$x(t) = A\cos\left(\int_0^{\omega t+\phi} sleaf_2(u)du\right) + A\sin\left(\int_0^{\omega t+\phi} sleaf_2(u)du\right) \tag{III.1}$$

With the above equation, the first order differential is obtained as follows:

$$\begin{aligned}\frac{dx(t)}{dt} &= -A\sin\left(\int_0^{\omega t+\phi} sleaf_2(u)du\right)\cdot\omega\cdot sleaf_2(\omega t+\phi) \\ &+ A\cos\left(\int_0^{\omega t+\phi} sleaf_2(u)du\right)\cdot\omega\cdot sleaf_2(\omega t+\phi)\end{aligned} \tag{III.2}$$

With the above equation, the second order differential is obtained as follows:

$$\begin{aligned}\frac{d^2x(t)}{dt^2} &= -A\cos\left(\int_0^{\omega t+\phi} sleaf_2(u)du\right)\cdot\omega^2\cdot(sleaf_2(\omega t+\phi))^2 \\ &- A\sin\left(\int_0^{\omega t+\phi} sleaf_2(u)du\right)\cdot\omega^2\cdot\sqrt{1-(sleaf_2(\omega t+\phi))^4} \\ &- A\sin\left(\int_0^{\omega t+\phi} sleaf_2(u)du\right)\cdot\omega^2\cdot(sleaf_2(\omega t+\phi))^2 \\ &+ A\cos\left(\int_0^{\omega t+\phi} sleaf_2(u)du\right)\cdot\omega^2\cdot\sqrt{1-(sleaf_2(\omega t+\phi))^4}\end{aligned} \tag{III.3}$$

By using Eq. (63) in Ref. [13], the following equation is obtained:

$$(sleaf_2(\omega t+\phi))^2 = \sin\left(2\int_0^{\omega t+\phi} sleaf_2(u)du\right) \tag{III.4}$$

The above equation is transformed as follows:

$$\sqrt{1-(sleaf_2(\omega t+\phi))^4} = \sqrt{1-\left(\sin\left(2\int_0^{\omega t+\phi} sleaf_2(u)du\right)\right)^2} = \cos\left(2\int_0^{\omega t+\phi} sleaf_2(u)du\right) \tag{III.5}$$

Substituting the above two equations into Eq. (III.3), the following equation is obtained:

$$\begin{aligned}\frac{d^2x(t)}{dt^2} &= -A\omega^2 \sin\left(\int_0^{\omega t+\phi} sleaf_2(u)du\right)\cdot\sin\left(2\int_0^{\omega t+\phi} sleaf_2(u)du\right) \\ &+ A\omega^2 \cos\left(\int_0^{\omega t+\phi} sleaf_2(u)du\right)\cdot\cos\left(2\int_0^{\omega t+\phi} sleaf_2(u)du\right) \\ &- A\omega^2 \cos\left(\int_0^{\omega t+\phi} sleaf_2(u)du\right)\cdot\sin\left(2\int_0^{\omega t+\phi} sleaf_2(u)du\right) \\ &- A\omega^2 \sin\left(\int_0^{\omega t+\phi} sleaf_2(u)du\right)\cdot\cos\left(2\int_0^{\omega t+\phi} sleaf_2(u)du\right)\end{aligned} \qquad (\text{III}.6)$$

The above equation is summarized by the addition theorem:

$$\begin{aligned}\frac{d^2x(t)}{dt^2} &= A\omega^2 \cos\left(2\int_0^{\omega t+\phi} sleaf_2(u)du + \int_0^{\omega t+\phi} sleaf_2(u)du\right) \\ &- A\omega^2 \sin\left(2\int_0^{\omega t+\phi} sleaf_2(u)du + \int_0^{\omega t+\phi} sleaf_2(u)du\right) \\ &= A\omega^2 \cos\left(3\int_0^{\omega t+\phi} sleaf_2(u)du\right) - A\omega^2 \sin\left(3\int_0^{\omega t+\phi} sleaf_2(u)du\right)\end{aligned} \qquad (\text{III}.7)$$

The above equation is transformed as follows:

$$\begin{aligned}\frac{d^2x(t)}{dt^2} &= 4A\omega^2\left\{\cos\left(\int_0^{\omega t+\phi} sleaf_2(u)du\right)\right\}^3 + 4A\omega^2\left\{\sin\left(\int_0^{\omega t+\phi} sleaf_2(u)du\right)\right\}^3 \\ &- 3A\omega^2 \cos\left(\int_0^{\omega t+\phi} sleaf_2(u)du\right) - 3A\omega^2 \sin\left(\int_0^{\omega t+\phi} sleaf_2(u)du\right)\end{aligned} \qquad (\text{III}.8)$$

Here, the following equation is derived.

$$\begin{aligned}\{x(t)\}^2 &= A^2\left(\sin\left(\int_0^{\omega t+\phi} sleaf_2(u)du\right)\right)^2 + A^2\left(\cos\left(\int_0^{\omega t+\phi} sleaf_2(u)du\right)\right)^2 \\ &+ 2A^2 \sin\left(\int_0^{\omega t+\phi} sleaf_2(u)du\right)\cos\left(\int_0^{\omega t+\phi} sleaf_2(u)du\right) \\ &= A^2 + 2A^2 \sin\left(\int_0^{\omega t+\phi} sleaf_2(u)du\right)\cos\left(\int_0^{\omega t+\phi} sleaf_2(u)du\right)\end{aligned} \qquad (\text{III}.9)$$

The following equation is derived from the above equation:

$$\sin\left(\int_0^{\omega t+\phi} sleaf_2(u)du\right)\cos\left(\int_0^{\omega t+\phi} sleaf_2(u)du\right)=\frac{\{x(t)\}^2-A^2}{2A^2} \quad (\text{III}.10)$$

Next, the following equation is derived.

$$\begin{aligned}\{x(t)\}^3 &= \left\{A\sin\left(\int_0^{\omega t+\phi} sleaf_2(u)du\right)+A\cos\left(\int_0^{\omega t+\phi} sleaf_2(u)du\right)\right\}^3 \\ &= A^3\left(\sin\left(\int_0^{\omega t+\phi} sleaf_2(u)du\right)\right)^3 + A^3\cdot 3\left(\sin\left(\int_0^{\omega t+\phi} sleaf_2(u)du\right)\right)^2\cdot\cos\left(\int_0^{\omega t+\phi} sleaf_2(u)du\right) \\ &+ A^3\cdot 3\cdot\sin\left(\int_0^{\omega t+\phi} sleaf_2(u)du\right)\cdot\left(\cos\left(\int_0^{\omega t+\phi} sleaf_2(u)du\right)\right)^2 + A^3\left(\cos\left(\int_0^{\omega t+\phi} sleaf_2(u)du\right)\right)^3 \\ &= A^3\left(\sin\left(\int_0^{\omega t+\phi} sleaf_2(u)du\right)\right)^3 + A^3\left(\cos\left(\int_0^{\omega t+\phi} sleaf_2(u)du\right)\right)^3 \\ &+ A^2\cdot 3\sin\left(\int_0^{\omega t+\phi} sleaf_2(u)du\right)\cdot\cos\left(\int_0^{\omega t+\phi} sleaf_2(u)du\right) \\ &\cdot\left(A\sin\left(\int_0^{\omega t+\phi} sleaf_2(u)du\right)+A\cos\left(\int_0^{\omega t+\phi} sleaf_2(u)du\right)\right)\end{aligned}$$

$$(\text{III}.11)$$

The following equation is derived by transforming the above equation:

$$\begin{aligned}&A^3\left(\sin\left(\int_0^{\omega t+\phi} sleaf_2(u)du\right)\right)^3 + A^3\left(\cos\left(\int_0^{\omega t+\phi} sleaf_2(u)du\right)\right)^3 \\ &= \{x(t)\}^3 - A^2\cdot 3\sin\left(\int_0^{\omega t+\phi} sleaf_2(u)du\right)\cdot\cos\left(\int_0^{\omega t+\phi} sleaf_2(u)du\right) \\ &\cdot\left(A\sin\left(\int_0^{\omega t+\phi} sleaf_2(u)du\right)+A\cos\left(\int_0^{\omega t+\phi} sleaf_2(u)du\right)\right) \\ &= \{x(t)\}^3 - A^2\cdot 3\frac{\{x(t)\}^2-A^2}{2A^2}x(t)=-\frac{1}{2}\{x(t)\}^3+\frac{3}{2}A^2 x(t)\end{aligned}$$

$$(\text{III}.12)$$

Substituting Eq. (III.12) into Eq. (III.11), the following equation is derived:

$$\begin{aligned}
\frac{d^2x(t)}{dt^2} &= 4\frac{\omega^2}{A^2}\left[A^3\left\{\cos\left(\int_0^{\omega t+\phi} sleaf_2(u)du\right)\right\}^3 + A^3\left\{\sin\left(\int_0^{\omega t+\phi} sleaf_2(u)du\right)\right\}^3\right] \\
&\quad -3\omega^2\left[A\cos\left(\int_0^{\omega t+\phi} sleaf_2(u)du\right) + A\sin\left(\int_0^{\omega t+\phi} sleaf_2(u)du\right)\right] \\
&= 4\frac{\omega^2}{A^2}\left[-\frac{1}{2}\{x(t)\}^3 + \frac{3}{2}A^2x(t)\right] - 3\omega^2 x(t) \\
&= -2\frac{\omega^2}{A^2}\{x(t)\}^3 + 3\omega^2 x(t)
\end{aligned} \qquad (\text{III}.13)$$

Eq. (2.11) is obtained from the above equation.

$$\frac{d^2x(t)}{dt^2} - 3\omega^2 x(t) + 2\left(\frac{\omega}{A}\right)^2 x(t)^3 = 0 \qquad (\text{III}.14)$$

The coefficients $\alpha$ and $\beta$ in Eq. (2.1) are expressed as follows:

$$\alpha = -3\omega^2 \qquad (\text{III}.15)$$

$$\beta = 2\left(\frac{\omega}{A}\right)^2 \qquad (\text{III}.16)$$

The initial condition of Eq. (2.12) and Eq. (2.13) is given as follows:

$$x(0) = A\cos\left(\int_0^\phi sleaf_2(u)du\right) + A\sin\left(\int_0^\phi sleaf_2(u)du\right) = \sqrt{2}A\cos\left(\int_0^\phi sleaf_2(u)du - \frac{\pi}{4}\right) \qquad (\text{III}.17)$$

$$\begin{aligned}
\frac{dx(0)}{dt} &= -A\sin\left(\int_0^\phi sleaf_2(u)du\right)\cdot\omega\cdot sleaf_2(\phi) + A\cos\left(\int_0^\phi sleaf_2(u)du\right)\cdot\omega\cdot sleaf_2(\phi) \\
&= A\cdot\omega\cdot sleaf_2(\phi)\cdot\left\{\cos\left(\int_0^\phi sleaf_2(u)du\right) - \sin\left(\int_0^\phi sleaf_2(u)du\right)\right\} \\
&= \sqrt{2}A\cdot\omega\cdot sleaf_2(\phi)\cdot\cos\left(\int_0^\phi sleaf_2(u)du + \frac{\pi}{4}\right)
\end{aligned} \qquad (\text{III}.18)$$

*Appendix IV*

The type(IV) exact solution is satisfied with the cubic duffing equation in Eq. (2.1). Eq. (2.14) is given as follows:

$$x(t) = A\cos\left(\int_0^{\omega t+\phi} sleaf_2(u)du\right) - A\sin\left(\int_0^{\omega t+\phi} sleaf_2(u)du\right) \tag{IV.1}$$

With the above equation, the first order differential is obtained as follows:

$$\begin{aligned}\frac{dx(t)}{dt} &= -A\sin\left(\int_0^{\omega t+\phi} sleaf_2(u)du\right)\cdot\omega\cdot sleaf_2(\omega t+\phi) \\ &\quad - A\cos\left(\int_0^{\omega t+\phi} sleaf_2(u)du\right)\cdot\omega\cdot sleaf_2(\omega t+\phi)\end{aligned} \tag{IV.2}$$

With the above equation, the second order differential is obtained as follows:

$$\begin{aligned}\frac{d^2x(t)}{dt^2} &= -A\cos\left(\int_0^{\omega t+\phi} sleaf_2(u)du\right)\cdot\omega^2\cdot(sleaf_2(\omega t+\phi))^2 \\ &\quad - A\sin\left(\int_0^{\omega t+\phi} sleaf_2(u)du\right)\cdot\omega^2\cdot\sqrt{1-(sleaf_2(\omega t+\phi))^4} \\ &\quad + A\sin\left(\int_0^{\omega t+\phi} sleaf_2(u)du\right)\cdot\omega^2\cdot(sleaf_2(\omega t+\phi))^2 \\ &\quad - A\cos\left(\int_0^{\omega t+\phi} sleaf_2(u)du\right)\cdot\omega^2\cdot\sqrt{1-(sleaf_2(\omega t+\phi))^4}\end{aligned} \tag{IV.3}$$

Substituting Eqs. (III.4) and (III.5) into the above equations, the following equation is obtained:

$$\begin{aligned}\frac{d^2x(t)}{dt^2} &= -A\omega^2\cos\left(\int_0^{\omega t+\phi} sleaf_2(u)du\right)\cdot\sin\left(2\int_0^{\omega t+\phi} sleaf_2(u)du\right) \\ &\quad - A\omega^2\sin\left(\int_0^{\omega t+\phi} sleaf_2(u)du\right)\cdot\cos\left(2\int_0^{\omega t+\phi} sleaf_2(u)du\right) \\ &\quad + A\omega^2\sin\left(\int_0^{\omega t+\phi} sleaf_2(u)du\right)\cdot\sin\left(2\int_0^{\omega t+\phi} sleaf_2(u)du\right) \\ &\quad - A\omega^2\cos\left(\int_0^{\omega t+\phi} sleaf_2(u)du\right)\cdot\cos\left(2\int_0^{\omega t+\phi} sleaf_2(u)du\right)\end{aligned} \tag{IV.4}$$

The above equation is summarized by the addition theorem:

$$\begin{aligned}\frac{d^2x(t)}{dt^2} &= -A\omega^2 \cos\left(2\int_0^{\omega t+\phi} sleaf_2(u)du + \int_0^{\omega t+\phi} sleaf_2(u)du\right) \\ &\quad - A\omega^2 \sin\left(2\int_0^{\omega t+\phi} sleaf_2(u)du + \int_0^{\omega t+\phi} sleaf_2(u)du\right) \\ &= -A\omega^2 \cos\left(3\int_0^{\omega t+\phi} sleaf_2(u)du\right) - A\omega^2 \sin\left(3\int_0^{\omega t+\phi} sleaf_2(u)du\right)\end{aligned} \quad (\text{IV}.5)$$

The above equation is transformed as follows:

$$\begin{aligned}\frac{d^2x(t)}{dt^2} &= -4A\omega^2\left\{\cos\left(\int_0^{\omega t+\phi} sleaf_2(u)du\right)\right\}^3 + 4A\omega^2\left\{\sin\left(\int_0^{\omega t+\phi} sleaf_2(u)du\right)\right\}^3 \\ &\quad + 3A\omega^2 \cos\left(\int_0^{\omega t+\phi} sleaf_2(u)du\right) - 3A\omega^2 \sin\left(\int_0^{\omega t+\phi} sleaf_2(u)du\right)\end{aligned} \quad (\text{IV}.6)$$

Here, the following equation is derived.

$$\begin{aligned}\{x(t)\}^2 &= A^2\left(\cos\left(\int_0^{\omega t+\phi} sleaf_2(u)du\right)\right)^2 + A^2\left(\sin\left(\int_0^{\omega t+\phi} sleaf_2(u)du\right)\right)^2 \\ &\quad - 2A^2 \sin\left(\int_0^{\omega t+\phi} sleaf_2(u)du\right)\cos\left(\int_0^{\omega t+\phi} sleaf_2(u)du\right) \\ &= A^2 - 2A^2 \sin\left(\int_0^{\omega t+\phi} sleaf_2(u)du\right)\cos\left(\int_0^{\omega t+\phi} sleaf_2(u)du\right)\end{aligned} \quad (\text{IV}.7)$$

The following equation is derived from the above equation:

$$\sin\left(\int_0^{\omega t+\phi} sleaf_2(u)du\right)\cos\left(\int_0^{\omega t+\phi} sleaf_2(u)du\right) = \frac{A^2 - \{x(t)\}^2}{2A^2} \quad (\text{IV}.8)$$

Next, the following equation is derived.

$$\{x(t)\}^3 = \left\{A\cos\left(\int_0^{\omega t+\phi} sleaf_2(u)du\right) - A\sin\left(\int_0^{\omega t+\phi} sleaf_2(u)du\right)\right\}^3$$

$$= A^3\left(\cos\left(\int_0^{\omega t+\phi} sleaf_2(u)du\right)\right)^3 - A^3 \cdot 3 \cdot \sin\left(\int_0^{\omega t+\phi} sleaf_2(u)du\right) \cdot \left(\cos\left(\int_0^{\omega t+\phi} sleaf_2(u)du\right)\right)^2$$

$$+ A^3 \cdot 3\left(\sin\left(\int_0^{\omega t+\phi} sleaf_2(u)du\right)\right)^2 \cdot \cos\left(\int_0^{\omega t+\phi} sleaf_2(u)du\right) - A^3\left(\sin\left(\int_0^{\omega t+\phi} sleaf_2(u)du\right)\right)^3$$

$$= A^3\left(\cos\left(\int_0^{\omega t+\phi} sleaf_2(u)du\right)\right)^3 - A^3\left(\sin\left(\int_0^{\omega t+\phi} sleaf_2(u)du\right)\right)^3$$

$$- A^2 \cdot 3\sin\left(\int_0^{\omega t+\phi} sleaf_2(u)du\right) \cdot \cos\left(\int_0^{\omega t+\phi} sleaf_2(u)du\right)$$

$$\cdot \left(A\cos\left(\int_0^{\omega t+\phi} sleaf_2(u)du\right) - A\sin\left(\int_0^{\omega t+\phi} sleaf_2(u)du\right)\right)$$

(IV.9)

The following equation is derived by transforming the above equation:

$$A^3\left(\cos\left(\int_0^{\omega t+\phi} sleaf_2(u)du\right)\right)^3 - A^3\left(\sin\left(\int_0^{\omega t+\phi} sleaf_2(u)du\right)\right)^3$$

$$= \{x(t)\}^3 + A^2 \cdot 3\sin\left(\int_0^{\omega t+\phi} sleaf_2(u)du\right) \cdot \cos\left(\int_0^{\omega t+\phi} sleaf_2(u)du\right)$$

$$\cdot \left(A\cos\left(\int_0^{\omega t+\phi} sleaf_2(u)du\right) - A\sin\left(\int_0^{\omega t+\phi} sleaf_2(u)du\right)\right)$$

$$= \{x(t)\}^3 + A^2 \cdot 3\frac{A^2 - \{x(t)\}^2}{2A^2} x(t) = -\frac{1}{2}\{x(t)\}^3 + \frac{3}{2}A^2 x(t)$$

(IV.10)

The following equation is derived from Eq. (IV.6).

$$\frac{d^2 x(t)}{dt^2} = -4\frac{\omega^2}{A^2}\left[A^3\left\{\cos\left(\int_0^{\omega t+\phi} sleaf_2(u)du\right)\right\}^3 - A^3\left\{\sin\left(\int_0^{\omega t+\phi} sleaf_2(u)du\right)\right\}^3\right]$$

$$+ 3\omega^2\left[A\cos\left(\int_0^{\omega t+\phi} sleaf_2(u)du\right) - A\sin\left(\int_0^{\omega t+\phi} sleaf_2(u)du\right)\right]$$

$$= -4\frac{\omega^2}{A^2}\left[-\frac{1}{2}\{x(t)\}^3 + \frac{3}{2}A^2 x(t)\right] + 3\omega^2 x(t)$$

$$= 2\frac{\omega^2}{A^2}\{x(t)\}^3 - 3\omega^2 x(t)$$

(IV.11)

Eq. (2.15) is obtained from the above equation.

$$\frac{d^2x(t)}{dt^2} + 3\omega^2 x(t) - 2\left(\frac{\omega}{A}\right)^2 x(t)^3 = 0 \tag{IV.12}$$

The coefficients $\alpha$ and $\beta$ in Eq. (2.1) are expressed as follows:

$$\alpha = 3\omega^2 \tag{IV.13}$$

$$\beta = -2\left(\frac{\omega}{A}\right)^2 \tag{IV.14}$$

The initial condition of Eq. (2.16) and Eq. (2.17) is given as follows:

$$x(0) = A\cos\left(\int_0^\phi sleaf_2(u)du\right) - A\sin\left(\int_0^\phi sleaf_2(u)du\right) = \sqrt{2}A\cos\left(\int_0^\phi sleaf_2(u)du + \frac{\pi}{4}\right) \tag{IV.15}$$

$$\begin{aligned}\frac{dx(0)}{dt} &= -A\sin\left(\int_0^\phi sleaf_2(u)du\right) \cdot \omega \cdot sleaf_2(\phi) + A\cos\left(\int_0^\phi sleaf_2(u)du\right) \cdot \omega \cdot sleaf_2(\phi) \\ &= A \cdot \omega \cdot sleaf_2(\phi) \cdot \left\{\cos\left(\int_0^\phi sleaf_2(u)du\right) - \sin\left(\int_0^\phi sleaf_2(u)du\right)\right\} \\ &= \sqrt{2}A \cdot \omega \cdot sleaf_2(\phi) \cdot \cos\left(\int_0^\phi sleaf_2(u)du + \frac{\pi}{4}\right)\end{aligned} \tag{IV.16}$$

*Appendix V*

The type(V) exact solution is satisfied with the cubic duffing equation in Eq. (2.1). Eq. (2.19) is given as follows:

$$x(t) = A\cos\left(\int_0^{\omega t+\phi} sleaf_2(u)du\right) + A\sin\left(\int_0^{\omega t+\phi} sleaf_2(u)du\right) + \sqrt{2}A\cos\left(\int_0^{\omega t+\phi} cleaf_2(u)du\right) \tag{V.1}$$

With the above equation, the first order differential is obtained as follows:

$$\frac{dx(t)}{dt} = -A\sin\left(\int_0^{\omega t+\phi} sleaf_2(u)du\right)\cdot\omega\cdot sleaf_2(\omega t+\phi)$$
$$+ A\cos\left(\int_0^{\omega t+\phi} sleaf_2(u)du\right)\cdot\omega\cdot sleaf_2(\omega t+\phi)$$
$$- \sqrt{2}A\sin\left(\int_0^{\omega t+\phi} cleaf_2(u)du\right)\cdot\omega\cdot cleaf_2(\omega t+\phi)$$

(V.2)

With the above equation, the second order differential is obtained as follows:

$$\frac{d^2x(t)}{dt^2} = -A\cos\left(\int_0^{\omega t+\phi} sleaf_2(u)du\right)\cdot\omega^2\cdot(sleaf_2(\omega t+\phi))^2$$
$$- A\sin\left(\int_0^{\omega t+\phi} sleaf_2(u)du\right)\cdot\omega^2\cdot\sqrt{1-(sleaf_2(\omega t+\phi))^4}$$
$$- A\sin\left(\int_0^{\omega t+\phi} sleaf_2(u)du\right)\cdot\omega^2\cdot(sleaf_2(\omega t+\phi))^2$$
$$+ A\cos\left(\int_0^{\omega t+\phi} sleaf_2(u)du\right)\cdot\omega^2\cdot\sqrt{1-(sleaf_2(\omega t+\phi))^4}$$
$$- \sqrt{2}A\cos\left(\int_0^{\omega t+\phi} cleaf_2(u)du\right)\cdot\omega^2\cdot(cleaf_2(\omega t+\phi))^2$$
$$+ \sqrt{2}A\sin\left(\int_0^{\omega t+\phi} cleaf_2(u)du\right)\cdot\omega^2\cdot\sqrt{1-(cleaf_2(\omega t+\phi))^4}$$

(V.3)

The above equation can be transformed as follows:

$$\frac{d^2x(t)}{dt^2} = -A\omega^2\cos\left(\int_0^{\omega t+\phi} sleaf_2(u)du\right)\cdot\sin\left(2\int_0^{\omega t+\phi} sleaf_2(u)du\right)$$
$$- A\omega^2\sin\left(\int_0^{\omega t+\phi} sleaf_2(u)du\right)\cdot\cos\left(2\int_0^{\omega t+\phi} sleaf_2(u)du\right)$$
$$- A\omega^2\sin\left(\int_0^{\omega t+\phi} sleaf_2(u)du\right)\cdot\sin\left(2\int_0^{\omega t+\phi} sleaf_2(u)du\right)$$
$$+ A\omega^2\cos\left(\int_0^{\omega t+\phi} sleaf_2(u)du\right)\cdot\cos\left(2\int_0^{\omega t+\phi} sleaf_2(u)du\right)$$
$$- \sqrt{2}A\omega^2\cos\left(\int_0^{\omega t+\phi} cleaf_2(u)du\right)\cdot\cos\left(2\int_0^{\omega t+\phi} cleaf_2(u)du\right)$$
$$+ \sqrt{2}A\omega^2\sin\left(\int_0^{\omega t+\phi} cleaf_2(u)du\right)\cdot\sin\left(2\int_0^{\omega t+\phi} cleaf_2(u)du\right)$$

(V.4)

The above equation is summarized by the addition theorem:

$$\frac{d^2x(t)}{dt^2} = A\omega^2 \cos\left(2\int_0^{\omega t+\phi} sleaf_2(u)du + \int_0^{\omega t+\phi} sleaf_2(u)du\right)$$

$$- A\omega^2 \sin\left(2\int_0^{\omega t+\phi} sleaf_2(u)du + \int_0^{\omega t+\phi} sleaf_2(u)du\right)$$

$$- \sqrt{2}A\omega^2 \cos\left(\int_0^{\omega t+\phi} cleaf_2(u)du + 2\int_0^{\omega t+\phi} cleaf_2(u)du\right)$$

$$= A\omega^2 \cos\left(3\int_0^{\omega t+\phi} sleaf_2(u)du\right) - A\omega^2 \sin\left(3\int_0^{\omega t+\phi} sleaf_2(u)du\right) - \sqrt{2}A\omega^2 \cos\left(3\int_0^{\omega t+\phi} cleaf_2(u)du\right)$$

(V.5)

The above equation can be transformed as follows:

$$\frac{d^2x(t)}{dt^2} = 4\left(\frac{\omega}{A}\right)^2\left\{A\cos\left(\int_0^{\omega t+\phi} sleaf_2(u)du\right)\right\}^3 + 4\left(\frac{\omega}{A}\right)^2\left\{A\sin\left(\int_0^{\omega t+\phi} sleaf_2(u)du\right)\right\}^3$$

$$- 4\sqrt{2}\left(\frac{\omega}{A}\right)^2\left\{A\cos\left(\int_0^{\omega t+\phi} cleaf_2(u)du\right)\right\}^3 - 3A\omega^2 \cos\left(\int_0^{\omega t+\phi} sleaf_2(u)du\right)$$

$$- 3A\omega^2 \sin\left(\int_0^{\omega t+\phi} sleaf_2(u)du\right) + 3\sqrt{2}\omega^2 A\cos\left(\int_0^{\omega t+\phi} cleaf_2(u)du\right)$$

(V.6)

Here, the following equation is derived from Eq. (66) in Ref. [13].

$$(sleaf_2(\omega t+\phi))^2 + (cleaf_2(\omega t+\phi))^2 + (sleaf_2(\omega t+\phi))^2(cleaf_2(\omega t+\phi))^2 = 1$$

(V.7)

The following equation is obtained from the above equation:

$$\cos\left(\int_0^{\omega t+\phi} cleaf_2(u)du\right)\left\{\sin\left(\int_0^{\omega t+\phi} sleaf_2(u)du\right) + \cos\left(\int_0^{\omega t+\phi} sleaf_2(u)du\right)\right\} = 1$$

(V.8)

Here, the following is expanded:

$$\left\{ A\sin\left(\int_0^{\omega t+\phi} sleaf_2(u)du\right) + A\cos\left(\int_0^{\omega t+\phi} sleaf_2(u)du\right)\right\}^3$$

$$= A^3\left(\sin\left(\int_0^{\omega t+\phi} sleaf_2(u)du\right)\right)^3 + A^3 \cdot 3\left(\sin\left(\int_0^{\omega t+\phi} sleaf_2(u)du\right)\right)^2 \cdot \cos\left(\int_0^{\omega t+\phi} sleaf_2(u)du\right)$$

$$+ A^3 \cdot 3 \cdot \sin\left(\int_0^{\omega t+\phi} sleaf_2(u)du\right) \cdot \left(\cos\left(\int_0^{\omega t+\phi} sleaf_2(u)du\right)\right)^2 + A^3\left(\cos\left(\int_0^{\omega t+\phi} sleaf_2(u)du\right)\right)^3$$

$$= A^3\left(\sin\left(\int_0^{\omega t+\phi} sleaf_2(u)du\right)\right)^3 + A^3\left(\cos\left(\int_0^{\omega t+\phi} sleaf_2(u)du\right)\right)^3$$

$$+ 3A^3 \cos\left(\int_0^{\omega t+\phi} sleaf_2(u)du\right) - 3A^3\left\{\cos\left(\int_0^{\omega t+\phi} sleaf_2(u)du\right)\right\}^3$$

$$+ 3A^3 \sin\left(\int_0^{\omega t+\phi} sleaf_2(u)du\right) - 3A^3\left\{\sin\left(\int_0^{\omega t+\phi} sleaf_2(u)du\right)\right\}^3$$

(V.9)

Next, $x(t)^3$ is expanded as follows:

$$\{x(t)\}^3 = A^3\left\{\cos\left(\int_0^{\omega t+\phi} sleaf_2(u)du\right) + \sin\left(\int_0^{\omega t+\phi} sleaf_2(u)du\right)\right\}^3 + 2\sqrt{2}A^3\left(\cos\left(\int_0^{\omega t+\phi} cleaf_2(u)du\right)\right)^3$$

$$+ 3\sqrt{2}A^3 \cos\left(\int_0^{\omega t+\phi} cleaf_2(u)du\right)\left\{\cos\left(\int_0^{\omega t+\phi} sleaf_2(u)du\right) + \sin\left(\int_0^{\omega t+\phi} sleaf_2(u)du\right)\right\}^2$$

$$+ 3\left(\sqrt{2}\right)^2 A^3\left\{\cos\left(\int_0^{\omega t+\phi} cleaf_2(u)du\right)\right\}^2\left\{\cos\left(\int_0^{\omega t+\phi} sleaf_2(u)du\right) + \sin\left(\int_0^{\omega t+\phi} sleaf_2(u)du\right)\right\}$$

(V.10)

Substituting Eq. (V.8) into the above equations, the following equation is obtained:

$$\{x(t)\}^3 = A^3\left(\sin\left(\int_0^{\omega t+\phi} sleaf_2(u)du\right)\right)^3 + A^3\left(\cos\left(\int_0^{\omega t+\phi} sleaf_2(u)du\right)\right)^3 + 3A^3 \cos\left(\int_0^{\omega t+\phi} sleaf_2(u)du\right)$$

$$- 3A^3\left\{\cos\left(\int_0^{\omega t+\phi} sleaf_2(u)du\right)\right\}^3 + 3A^3 \sin\left(\int_0^{\omega t+\phi} sleaf_2(u)du\right) - 3A^3\left\{\sin\left(\int_0^{\omega t+\phi} sleaf_2(u)du\right)\right\}^3$$

$$+ 2\sqrt{2}A^3\left(\cos\left(\int_0^{\omega t+\phi} cleaf_2(u)du\right)\right)^3 + 3\sqrt{2}A^3 \cos\left(\int_0^{\omega t+\phi} sleaf_2(u)du\right) + 3\sqrt{2}A^3 \sin\left(\int_0^{\omega t+\phi} sleaf_2(u)du\right)$$

$$+ 6A^3 \cos\left(\int_0^{\omega t+\phi} cleaf_2(u)du\right)$$

(V.11)

The following equation is obtained from the above equation:

$$\{x(t)\}^3 = -2A^3\left(\sin\left(\int_0^{\omega t+\phi} sleaf_2(u)du\right)\right)^3 - 2A^3\left(\cos\left(\int_0^{\omega t+\phi} sleaf_2(u)du\right)\right)^3 + 2\sqrt{2}A^3\left(\cos\left(\int_0^{\omega t+\phi} cleaf_2(u)du\right)\right)^3$$
$$+ 3A^3\cos\left(\int_0^{\omega t+\phi} sleaf_2(u)du\right) + 3A^3\sin\left(\int_0^{\omega t+\phi} sleaf_2(u)du\right) + 3\sqrt{2}A^3\cos\left(\int_0^{\omega t+\phi} sleaf_2(u)du\right)$$
$$+ 3\sqrt{2}A^3\sin\left(\int_0^{\omega t+\phi} sleaf_2(u)du\right) + 6A^3\cos\left(\int_0^{\omega t+\phi} cleaf_2(u)du\right)$$

(V.12)

Eq. (2.19) is obtained from the above equation.

$$\frac{d^2x(t)}{dt^2} - 3\omega^2\left(1+2\sqrt{2}\right)x(t) + 2\frac{\omega^2}{A^2}x(t)^3 = 0 \tag{V.13}$$

The initial condition of Eq. (2.20) and Eq. (2.21) is given as follows:

$$x(0) = A\cos\left(\int_0^{\phi} sleaf_2(u)du\right) + A\sin\left(\int_0^{\phi} sleaf_2(u)du\right) + \sqrt{2}A\cos\left(\int_0^{\phi} cleaf_2(u)du\right)$$
$$= \sqrt{2}A\left\{\sin\left(\int_0^{\phi} sleaf_2(u)du + \frac{\pi}{4}\right) + \cos\left(\int_0^{\phi} cleaf_2(u)du\right)\right\} \tag{V.14}$$

$$\frac{dx(0)}{dt} = -A\sin\left(\int_0^{\phi} sleaf_2(u)du\right)\cdot\omega\cdot sleaf_2(\phi)$$
$$+ A\cos\left(\int_0^{\phi} sleaf_2(u)du\right)\cdot\omega\cdot sleaf_2(\phi)$$
$$- \sqrt{2}A\sin\left(\int_0^{\phi} cleaf_2(u)du\right)\cdot\omega\cdot cleaf_2(\phi)$$
$$= \sqrt{2}A\cdot\omega\cdot sleaf_2(\phi)\cdot\cos\left(\int_0^{\phi} sleaf_2(u)du + \frac{\pi}{4}\right) - \sqrt{2}A\sin\left(\int_0^{\phi} cleaf_2(u)du\right)\cdot\omega\cdot cleaf_2(\phi)$$
$$= \sqrt{2}A\cdot\omega\cdot\left\{\cos\left(\int_0^{\phi} sleaf_2(u)du + \frac{\pi}{4}\right)\cdot sleaf_2(\phi) - \sin\left(\int_0^{\phi} cleaf_2(u)du\right)\cdot cleaf_2(\phi)\right\}$$

(V.15)

*Appendix VI*

The type(VI) exact solution is satisfied with the cubic duffing equation in Eq. (2.1). Eq. (2.23) is given as follows:

$$x(t) = A\cos\left(\int_0^{\omega t+\phi} sleaf_2(u)du\right) + A\sin\left(\int_0^{\omega t+\phi} sleaf_2(u)du\right) - \sqrt{2}A\cos\left(\int_0^{\omega t+\phi} cleaf_2(u)du\right)$$

(VI.1)

With the above equation, the first order differential is obtained as follows:

$$\frac{dx(t)}{dt} = -A\sin\left(\int_0^{\omega t+\phi} sleaf_2(u)du\right) \cdot \omega \cdot sleaf_2(\omega t + \phi)$$
$$+ A\cos\left(\int_0^{\omega t+\phi} sleaf_2(u)du\right) \cdot \omega \cdot sleaf_2(\omega t + \phi)$$
$$+ \sqrt{2}A\sin\left(\int_0^{\omega t+\phi} cleaf_2(u)du\right) \cdot \omega \cdot cleaf_2(\omega t + \phi)$$

(VI.2)

With the above equation, the second order differential is obtained as follows:

$$\frac{d^2x(t)}{dt^2} = -A\cos\left(\int_0^{\omega t+\phi} sleaf_2(u)du\right) \cdot \omega^2 \cdot (sleaf_2(\omega t + \phi))^2$$
$$- A\sin\left(\int_0^{\omega t+\phi} sleaf_2(u)du\right) \cdot \omega^2 \cdot \sqrt{1 - (sleaf_2(\omega t + \phi))^4}$$
$$- A\sin\left(\int_0^{\omega t+\phi} sleaf_2(u)du\right) \cdot \omega^2 \cdot (sleaf_2(\omega t + \phi))^2$$
$$+ A\cos\left(\int_0^{\omega t+\phi} sleaf_2(u)du\right) \cdot \omega^2 \cdot \sqrt{1 - (sleaf_2(\omega t + \phi))^4}$$
$$+ \sqrt{2}A\cos\left(\int_0^{\omega t+\phi} cleaf_2(u)du\right) \cdot \omega^2 \cdot (cleaf_2(\omega t + \phi))^2$$
$$- \sqrt{2}A\sin\left(\int_0^{\omega t+\phi} cleaf_2(u)du\right) \cdot \omega^2 \cdot \sqrt{1 - (cleaf_2(\omega t + \phi))^4}$$

(VI.3)

The following equation is derived from the above equation:

$$\frac{d^2x(t)}{dt^2} = -A\omega^2 \cos\left(\int_0^{\omega t+\phi} sleaf_2(u)du\right)\cdot \sin\left(2\int_0^{\omega t+\phi} sleaf_2(u)du\right)$$

$$- A\omega^2 \sin\left(\int_0^{\omega t+\phi} sleaf_2(u)du\right)\cdot \cos\left(2\int_0^{\omega t+\phi} sleaf_2(u)du\right)$$

$$- A\omega^2 \sin\left(\int_0^{\omega t+\phi} sleaf_2(u)du\right)\cdot \sin\left(2\int_0^{\omega t+\phi} sleaf_2(u)du\right)$$

$$+ A\omega^2 \cos\left(\int_0^{\omega t+\phi} sleaf_2(u)du\right)\cdot \cos\left(2\int_0^{\omega t+\phi} sleaf_2(u)du\right)$$

$$+ \sqrt{2}A\omega^2 \cos\left(\int_0^{\omega t+\phi} cleaf_2(u)du\right)\cdot \cos\left(2\int_0^{\omega t+\phi} cleaf_2(u)du\right)$$

$$- \sqrt{2}A\omega^2 \sin\left(\int_0^{\omega t+\phi} cleaf_2(u)du\right)\cdot \sin\left(2\int_0^{\omega t+\phi} cleaf_2(u)du\right)$$

(VI.4)

The above equation is summarized by the addition theorem:

$$\frac{d^2x(t)}{dt^2} = A\omega^2 \cos\left(2\int_0^{\omega t+\phi} sleaf_2(u)du + \int_0^{\omega t+\phi} sleaf_2(u)du\right)$$

$$- A\omega^2 \sin\left(2\int_0^{\omega t+\phi} sleaf_2(u)du + \int_0^{\omega t+\phi} sleaf_2(u)du\right)$$

$$+ \sqrt{2}A\omega^2 \cos\left(\int_0^{\omega t+\phi} cleaf_2(u)du + 2\int_0^{\omega t+\phi} cleaf_2(u)du\right)$$

$$= A\omega^2 \cos\left(3\int_0^{\omega t+\phi} sleaf_2(u)du\right) - A\omega^2 \sin\left(3\int_0^{\omega t+\phi} sleaf_2(u)du\right) + \sqrt{2}A\omega^2 \cos\left(3\int_0^{\omega t+\phi} cleaf_2(u)du\right)$$

(VI.5)

The above equation can be transformed as follows:

$$\frac{d^2x(t)}{dt^2} = 4\left(\frac{\omega}{A}\right)^2 \left\{A\cos\left(\int_0^{\omega t+\phi} sleaf_2(u)du\right)\right\}^3 + 4\left(\frac{\omega}{A}\right)^2 \left\{A\sin\left(\int_0^{\omega t+\phi} sleaf_2(u)du\right)\right\}^3$$

$$+ 4\sqrt{2}\left(\frac{\omega}{A}\right)^2 \left\{A\cos\left(\int_0^{\omega t+\phi} cleaf_2(u)du\right)\right\}^3 - 3A\omega^2 \cos\left(\int_0^{\omega t+\phi} sleaf_2(u)du\right)$$

$$- 3A\omega^2 \sin\left(\int_0^{\omega t+\phi} sleaf_2(u)du\right) - 3\sqrt{2}\omega^2 A\cos\left(\int_0^{\omega t+\phi} cleaf_2(u)du\right)$$

(VI.6)

Next, $x(t)^3$ is expanded as follows:

$$\{x(t)\}^3 = A^3 \left\{ \cos\left(\int_0^{\omega t+\phi} sleaf_2(u)du\right) + \sin\left(\int_0^{\omega t+\phi} sleaf_2(u)du\right) \right\}^3 - 2\sqrt{2}A^3 \left( \cos\left(\int_0^{\omega t+\phi} cleaf_2(u)du\right) \right)^3$$

$$- 3\sqrt{2}A^3 \cos\left(\int_0^{\omega t+\phi} cleaf_2(u)du\right) \left\{ \cos\left(\int_0^{\omega t+\phi} sleaf_2(u)du\right) + \sin\left(\int_0^{\omega t+\phi} sleaf_2(u)du\right) \right\}^2$$

$$+ 3\left(\sqrt{2}\right)^2 A^3 \left\{ \cos\left(\int_0^{\omega t+\phi} cleaf_2(u)du\right) \right\}^2 \left\{ \cos\left(\int_0^{\omega t+\phi} sleaf_2(u)du\right) + \sin\left(\int_0^{\omega t+\phi} sleaf_2(u)du\right) \right\}$$

(VI.7)

Substituting Eq. (V.8) into the above equations, the following equation is obtained:

$$\{x(t)\}^3 = A^3 \left( \sin\left(\int_0^{\omega t+\phi} sleaf_2(u)du\right) \right)^3 + A^3 \left( \cos\left(\int_0^{\omega t+\phi} sleaf_2(u)du\right) \right)^3 + 3A^3 \cos\left(\int_0^{\omega t+\phi} sleaf_2(u)du\right)$$

$$- 3A^3 \left\{ \cos\left(\int_0^{\omega t+\phi} sleaf_2(u)du\right) \right\}^3 + 3A^3 \sin\left(\int_0^{\omega t+\phi} sleaf_2(u)du\right) - 3A^3 \left\{ \sin\left(\int_0^{\omega t+\phi} sleaf_2(u)du\right) \right\}^3$$

$$- 2\sqrt{2}A^3 \left( \cos\left(\int_0^{\omega t+\phi} cleaf_2(u)du\right) \right)^3 - 3\sqrt{2}A^3 \cos\left(\int_0^{\omega t+\phi} sleaf_2(u)du\right) - 3\sqrt{2}A^3 \sin\left(\int_0^{\omega t+\phi} sleaf_2(u)du\right)$$

$$+ 6A^3 \cos\left(\int_0^{\omega t+\phi} cleaf_2(u)du\right)$$

(VI.8)

The following equation is obtained from the above equation:

$$\{x(t)\}^3 = -2A^3 \left( \sin\left(\int_0^{\omega t+\phi} sleaf_2(u)du\right) \right)^3 - 2A^3 \left( \cos\left(\int_0^{\omega t+\phi} sleaf_2(u)du\right) \right)^3 - 2\sqrt{2}A^3 \left( \cos\left(\int_0^{\omega t+\phi} cleaf_2(u)du\right) \right)^3$$

$$+ 3A^3 \cos\left(\int_0^{\omega t+\phi} sleaf_2(u)du\right) + 3A^3 \sin\left(\int_0^{\omega t+\phi} sleaf_2(u)du\right) - 3\sqrt{2}A^3 \cos\left(\int_0^{\omega t+\phi} sleaf_2(u)du\right)$$

$$- 3\sqrt{2}A^3 \sin\left(\int_0^{\omega t+\phi} sleaf_2(u)du\right) + 6A^3 \cos\left(\int_0^{\omega t+\phi} cleaf_2(u)du\right)$$

(VI.9)

Eq. (2.23) is obtained from the above equation.

$$\frac{d^2 x(t)}{dt^2} + 3\omega^2 \left(2\sqrt{2}-1\right) x(t) + 2\frac{\omega^2}{A^2} x(t)^3 = 0$$

(VI.10)

The initial condition of Eq. (2.24) and Eq. (2.25) is expressed as follows:

$$x(0) = A\cos\left(\int_0^\phi sleaf_2(u)du\right) + A\sin\left(\int_0^\phi sleaf_2(u)du\right) - \sqrt{2}A\cos\left(\int_0^\phi cleaf_2(u)du\right)$$

$$= \sqrt{2}A\left\{\sin\left(\int_0^\phi sleaf_2(u)du + \frac{\pi}{4}\right) - \cos\left(\int_0^\phi cleaf_2(u)du\right)\right\}$$

(VI.11)

$$\frac{dx(0)}{dt} = -A\sin\left(\int_0^\phi sleaf_2(u)du\right) \cdot \omega \cdot sleaf_2(\phi)$$

$$+ A\cos\left(\int_0^\phi sleaf_2(u)du\right) \cdot \omega \cdot sleaf_2(\phi)$$

$$+ \sqrt{2}A\sin\left(\int_0^\phi cleaf_2(u)du\right) \cdot \omega \cdot cleaf_2(\phi)$$

$$= \sqrt{2}A \cdot \omega \cdot sleaf_2(\phi) \cdot \cos\left(\int_0^\phi sleaf_2(u)du + \frac{\pi}{4}\right) + \sqrt{2}A\sin\left(\int_0^\phi cleaf_2(u)du\right) \cdot \omega \cdot cleaf_2(\phi)$$

$$= \sqrt{2}A \cdot \omega \cdot \left\{\cos\left(\int_0^\phi sleaf_2(u)du + \frac{\pi}{4}\right) \cdot sleaf_2(\phi) + \sin\left(\int_0^\phi cleaf_2(u)du\right) \cdot cleaf_2(\phi)\right\}$$

(VI.12)

## *Appendix VII*

The type(VII) exact solution is satisfied with the cubic duffing equation in Eq. (2.1). Eq. (2.26) is given as follows:

$$x(t) = A \cdot sleaf_2(\omega \cdot t + \phi) \cdot cleaf_2(\omega \cdot t + \phi)$$

(VII.1)

With the above equation, the first order differential is obtained as follows:

$$\frac{dx(t)}{dt} = A\omega \cdot cleaf_2(\omega \cdot t + \phi) \cdot \sqrt{1 - (sleaf_2(\omega \cdot t + \phi))^4} - A\omega \cdot sleaf_2(\omega \cdot t + \phi)\sqrt{1 - (cleaf_2(\omega \cdot t + \phi))^4}$$

(VII.2)

With the above equation, the second order differential is obtained as follows:

$$\frac{d^2x(t)}{dt^2} = -A\omega^2 \cdot \sqrt{1-(cleaf_2(\omega \cdot t + \phi))^4} \cdot \sqrt{1-(sleaf_2(\omega \cdot t + \phi))^4}$$
$$+ A\omega^2 \cdot cleaf_2(\omega \cdot t + \phi) \cdot \left\{-2(sleaf_2(\omega \cdot t + \phi))^3\right\}$$
$$- A\omega^2 \cdot \sqrt{1-(sleaf_2(\omega \cdot t + \phi))^4} \sqrt{1-(cleaf_2(\omega \cdot t + \phi))^4}$$
$$- A\omega^2 \cdot sleaf_2(\omega \cdot t + \phi)\left\{2(cleaf_2(\omega \cdot t + \phi))^3\right\}$$
$$= -2A\omega^2 \cdot \left\{sleaf_2(\omega \cdot t + \phi) \cdot cleaf_2(\omega \cdot t + \phi)\left((cleaf_2(\omega \cdot t + \phi))^2 + (sleaf_2(\omega \cdot t + \phi))^2\right)\right.$$
$$\left. - 2A\omega^2 \cdot \sqrt{1-(cleaf_2(\omega \cdot t + \phi))^4} \cdot \sqrt{1-(sleaf_2(\omega \cdot t + \phi))^4}\right\}$$

(Ⅶ.3)

The following equation is obtained from Eq. (Ⅴ.7).

$$\sqrt{1-(sleaf_2(\omega \cdot t + \phi))^4} = \sqrt{1-\left\{\frac{1-(cleaf_2(\omega \cdot t + \phi))^2}{1+(cleaf_2(\omega \cdot t + \phi))^2}\right\}^2}$$
$$= \sqrt{\frac{\{1+(cleaf_2(\omega \cdot t + \phi))^2\}^2 - \{1-(cleaf_2(\omega \cdot t + \phi))^2\}^2}{\{1+(cleaf_2(\omega \cdot t + \phi))^2\}^2}}$$
$$= \sqrt{\frac{4(cleaf_2(\omega \cdot t + \phi))^2}{\{1+(cleaf_2(\omega \cdot t + \phi))^2\}^2}} = \frac{2cleaf_2(\omega \cdot t + \phi)}{1+(cleaf_2(\omega \cdot t + \phi))^2}$$

(Ⅶ.4)

Similarly, the following equation is derived.

$$\sqrt{1-(cleaf_2(\omega \cdot t + \phi))^4} = \frac{2sleaf_2(\omega \cdot t + \phi)}{1+(sleaf_2(\omega \cdot t + \phi))^2}$$

(Ⅶ.5)

The following equation is derived by using these above equations.

$$\sqrt{1-(cleaf_2(\omega \cdot t + \phi))^4} \cdot \sqrt{1-(sleaf_2(\omega \cdot t + \phi))^4} = \frac{2sleaf_2(\omega \cdot t + \phi)}{1+(sleaf_2(\omega \cdot t + \phi))^2} \frac{2cleaf_2(\omega \cdot t + \phi)}{1+(cleaf_2(\omega \cdot t + \phi))^2}$$
$$= \frac{4sleaf_2(\omega \cdot t + \phi) \cdot cleaf_2(\omega \cdot t + \phi)}{1+(cleaf_2(\omega \cdot t + \phi))^2 + (sleaf_2(\omega \cdot t + \phi))^2 + (cleaf_2(\omega \cdot t + \phi))^2(sleaf_2(\omega \cdot t + \phi))^2}$$
$$= 2sleaf_2(\omega \cdot t + \phi) \cdot cleaf_2(\omega \cdot t + \phi)$$

(Ⅶ.6)

The following equation is derived by using these above equations.

$$\begin{aligned}\frac{d^2x(t)}{dt^2} &= -2A\omega^2 \cdot sleaf_2(\omega \cdot t + \phi) \cdot cleaf_2(\omega \cdot t + \phi)\{1-(cleaf_2(\omega \cdot t + \phi))^2(sleaf_2(\omega \cdot t + \phi))^2\} \\ &\quad - 4A\omega^2 \cdot sleaf_2(\omega \cdot t + \phi) \cdot cleaf_2(\omega \cdot t + \phi) \\ &= -6A\omega^2 \cdot sleaf_2(\omega \cdot t + \phi) \cdot cleaf_2(\omega \cdot t + \phi) + 2A\omega^2 \cdot (cleaf_2(\omega \cdot t + \phi))^3(sleaf_2(\omega \cdot t + \phi))^3 \\ &= -6\omega^2 \cdot \{A \cdot sleaf_2(\omega \cdot t + \phi) \cdot cleaf_2(\omega \cdot t + \phi)\} + 2\frac{\omega^2}{A^2} \cdot \{A \cdot sleaf_2(\omega \cdot t + \phi) \cdot cleaf_2(\omega \cdot t + \phi)\}^3\end{aligned}$$

(Ⅶ.7)

The above equation is transformed as follows:

$$\frac{d^2x(t)}{dt^2} + 6\omega^2 x(t) - 2\left(\frac{\omega}{A}\right)^2 x(t)^3 = 0 \tag{Ⅶ.8}$$

The coefficients $\alpha$ and $\beta$ in Eq. (2.1) are expressed as follows:

$$\alpha = 6\omega^2 \tag{Ⅶ.9}$$

$$\beta = -2\left(\frac{\omega}{A}\right)^2 \tag{Ⅶ.10}$$

The initial condition of Eq. (2.28) and Eq. (2.29) is given as follows:

$$x(0) = A \cdot sleaf_2(\phi) \cdot cleaf_2(\phi) \tag{Ⅶ.11}$$

$$\begin{aligned}\frac{dx(0)}{dt} &= A\omega \cdot cleaf_2(\phi) \cdot \sqrt{1-(sleaf_2(\phi))^4} - A\omega \cdot sleaf_2(\phi)\sqrt{1-(cleaf_2(\phi))^4} \\ &= A\omega \cdot cleaf_2(\phi) \cdot \frac{2cleaf_2(\phi)}{1+(cleaf_2(\phi))^2} - A\omega \cdot sleaf_2(\phi)\frac{2sleaf_2(\phi)}{1+(sleaf_2(\phi))^2} \\ &= A\omega \cdot (cleaf_2(\phi))^2 \cdot (1+(sleaf_2(\phi))^2) - A\omega \cdot (sleaf_2(\phi))^2(1+(cleaf_2(\phi))^2) \\ &= A\omega\{(cleaf_2(\phi))^2 - (sleaf_2(\phi))^2\}\end{aligned}$$

(Ⅶ.12)

## Appendix Ⅷ

We discuss the following derivative with respect to the variable $t$.

$$\frac{d}{dt}\arccos(cleaf_n(t))^n = -\frac{1}{\sqrt{1-(cleaf_n(t))^{2n}}} n(cleaf_n(t))^{n-1}\left\{-\sqrt{1-(cleaf_n(t))^{2n}}\right\} = n(cleaf_n(t))^{n-1} \tag{Ⅷ.1}$$

The above equation is integrated from 0 to the variable $t$.

$$\left[\arccos(cleaf_n(t))^n\right]_0^l = \int_0^l n(cleaf_n(t))^{n-1} dt \tag{VIII.2}$$

$$\left[\arccos(cleaf_n(t))^n\right]_0^l = \arccos(cleaf_n(l))^n - \arccos(cleaf_n(0))^n$$
$$= \arccos(cleaf_n(l))^n - \arccos(1) = \arccos(cleaf_n(l))^n \tag{VIII.3}$$

Therefore, the following relation is obtained as follows:

$$(cleaf_n(t))^n = \cos\left(n\int_0^t (cleaf_n(u))^{n-1} du\right) \tag{VIII.4}$$

In case of $n=2$, the above equation is as follows:

$$(cleaf_2(t))^2 = \cos\left(2\int_0^t cleaf_2(u) du\right) \tag{VIII.5}$$

Substituting $t = \frac{\pi_2}{2}(2m-1)$ into the above equation as follows:

$$\cos\left(2\int_0^{\frac{\pi_2}{2}(2m-1)} cleaf_2(u) du\right) = \left(cleaf_2\left(\frac{\pi_2}{2}(2m-1)\right)\right)^2 = 0 \tag{VIII.6}$$

The Eq. (1.17) is applied to the above equation. As shown in Fig.1, the range of $\int_0^t cleaf_2(u) du$ is as follows:

$$-\frac{\pi}{4} \leq \int_0^t cleaf_2(u) du \leq \frac{\pi}{4} \tag{VIII.7}$$

Therefore, the following equation based on the relation (VIII.6) can be obtained:

$$2\int_0^{\frac{\pi_2}{2}(2m-1)} cleaf_2(u) du = \pm\frac{\pi}{2} \tag{VIII.8}$$

The Eq. (3.1) can be obtained. Next, substituting $t=m\pi_2$ into the Eq. (VIII.5) as follows:

$$\cos\left(2\int_0^{m\pi_2} cleaf_2(t)dt\right) = (cleaf_2(m\pi_2))^2 = (\pm 1)^2 = 1 \tag{VIII.9}$$

The Eq. (1.18) or the Eq. (1.19) are applied to the above equation. In order to satisfy with the inequality (VIII.7), the integral: $\int_0^{m\pi_2} cleaf_2(t)dt$ become zero.